\documentclass[letterpaper,11pt]{article}

\newcommand*{\affaddr}[1]{#1} 

\newcommand*{\email}[1]{\texttt{#1}}

\usepackage{xr-hyper}
\usepackage{natbib}
\usepackage{hyperref} 
\usepackage{amsmath}
\usepackage{amssymb} 
\usepackage{amsthm}
\usepackage{geometry} 
\usepackage[dvipsnames]{xcolor}  
\usepackage{todonotes}
\usepackage{algorithm}
\usepackage{algorithmic}
\usepackage{ulem}
\usepackage{cancel}
\usepackage{comment}
\usepackage{caption}
\usepackage{subcaption}
\usepackage{todonotes}

\def\btheta{\boldsymbol\bx}
\def\grad{\nabla}

\def\bs{\mathbf{s}}

\def\bu{\mathbf{u}}
\def\bv{\mathbf{v}}

\def\bx{\mathbf{x}}  
\def\by{\mathbf{y}}

\def\cA{\mathcal{A}}

\def\cE{\mathcal{E}}
\def\cF{\mathcal{F}}

\def\cI{\mathcal{I}}

\def\cM{\mathcal{M}}
\def\cN{\mathcal{N}}
\def\cO{\mathcal{O}}

\def\cT{\mathcal{T}}
\def\cU{\mathcal{U}}
\def\cV{\mathcal{V}}

\def\mE{\mathbb{E}}

\def\mR{\mathbb{R}}

\def\smskip{\smallskip}

\def\texitem#1{\par\smskip\noindent\hangindent 25pt
               \hbox to 25pt {\hss #1 ~}\ignorespaces}

\def\norm#1{\|#1\|}

\newcommand{\BEAS}{\begin{eqnarray*}}
\newcommand{\EEAS}{\end{eqnarray*}}
\newcommand{\BEA}{\begin{eqnarray}}
\newcommand{\EEA}{\end{eqnarray}}
\newcommand{\BEQ}{\begin{eqnarray}}
\newcommand{\EEQ}{\end{eqnarray}}
\newcommand{\BIT}{\begin{itemize}}
\newcommand{\EIT}{\end{itemize}}
\newcommand{\BNUM}{\begin{enumerate}}
\newcommand{\ENUM}{\end{enumerate}}

\newcommand{\BA}{\begin{array}}
\newcommand{\EA}{\end{array}}

\newif\ifpagenumbering
\pagenumberingtrue

\pagenumberingfalse

\def\fprod#1{\left\langle#1\right\rangle}

\newtheorem{theorem}{Theorem}[section]

\newtheorem{lemma}{Lemma}[section]
\newtheorem{remark}{Remark}[section]

\newtheorem{definition}{Definition}
\newtheorem{assumption}{Assumption}[section]

\usepackage{xcolor}
\usepackage{enumerate}

\usepackage{relsize}
\usepackage{subfiles}

\usepackage{resizegather} 
\allowdisplaybreaks
\usepackage{adjustbox}

\newcounter{relctr} 
\everydisplay\expandafter{\the\everydisplay\setcounter{relctr}{0}} 

\newcommand\labelrel[2]{%
  \begingroup
    \refstepcounter{relctr}%
    \stackrel{\textnormal{(\alph{relctr})}}{\mathstrut{#1}}%
    \originallabel{#2}%
  \endgroup
}
\AtBeginDocument{\let\originallabel\label} 

\usepackage{natbib}
 \bibpunct[, ]{(}{)}{,}{a}{}{,}%
\begin{document}
\date{\today}
\title{\vspace{-1cm} \Large\bfseries Riemannian Stochastic Gradient Method for Nested Composition Optimization}

%

\author{%
	Dewei Zhang and Sam Davanloo Tajbakhsh \\
	\email{\{zhang.8705,davanloo.1\}@osu.edu}\\
	\vspace{0.5cm}
	\affaddr{The Ohio State University}%
}

\maketitle


\begin{abstract}
This work considers optimization of composition of functions in a nested form over Riemannian manifolds where each function contains an expectation.
This type of problems is gaining popularity in applications such as policy evaluation in reinforcement learning or model customization in meta-learning. The standard Riemannian stochastic gradient methods for non-compositional optimization cannot be directly applied as the stochastic approximation of the inner functions create biases in the gradients of the outer functions. For two-level composition optimization, we present a Riemannian Stochastic Composition Gradient Descent (R-SCGD) method that finds an approximate stationary point, with expected squared Riemannian gradient smaller than $\epsilon$, in $\cO(\epsilon^{-2})$ calls to the stochastic gradient oracle of the outer function and stochastic function and gradient oracles of the inner function. Furthermore, we generalize the R-SCGD algorithms for problems with multi-level nested compositional structures, with the same complexity of $\cO(\epsilon^{-2})$ for the first-order stochastic oracle. Finally, the performance of the R-SCGD method is numerically evaluated over a policy evaluation problem in reinforcement learning.
\end{abstract}

\section{Introduction}\label{sec:intro}
\vspace{-0.2cm}
We consider optimizing the nested stochastic composition of functions over Riemannian manifolds. The two-level composition optimization problem has the form 
\begin{equation}\label{eq:main1}
\min_{\bx\in \cM} F(\bx)= \mE_\xi[f_{\xi}(\mE_\phi[g_{\phi}(\bx)])]
\end{equation}
where $g_{\phi}:\cM\to\cN$ is a smooth map, $\cM$ and $\cN\subseteq\cE$ are Riemannian manifolds,  $\cE$ is a Euclidean space, $f_{\xi}:\cE\to\mR$ is a continuously differentiable function, and $\xi$ and $\phi$ are independent random variables \citep{absil2009optimization}. With a slight abuse of notation, we denote $\mE_{\xi}[f_{\xi}(\cdot)]$, $\mE_{\phi}[g_{\phi}(\cdot)]$ by $f(\cdot)$ and $g(\cdot)$, respectively. We assume throughout the paper that there exists at least one global optimal solution $\bx^*\in\cM$ to problem~\eqref{eq:main1}. We do not require either the outer function $f_{\xi}$ or the inner function $g_{\phi}$ to be (geodesically) convex or monotone. As a result, the composition problem cannot be reformulated into a saddle point problem in general~\citep{zhang2020optimal}. The Riemannian stochastic gradient method \citep{bonnabel2013stochastic} cannot be applied to solve \eqref{eq:main1} since the stochastic gradient $\nabla f_\xi$ evaluated at the stochastic function value $g_\phi$ does not result in an unbiased estimate of $\grad f(g(x))$ which then results in a biased estimate of $\nabla F(\bx)$~\citep{robbins1951stochastic}. Hence, a special algorithmic setup is required to tackle this type of problem.  
 
The stochastic compositional optimization finds applications in reinforcement learning for policy evaluation, meta-learning, stochastic minimax problems, dynamic programming and risk-averse problems~\citep{ben2002robust,bertsekas2012dynamic,goodfellow2016deep,finn2017model,sutton2018reinforcement}. The Riemannian manifold constraint arises due to 1) recent reformulation of some nonconvex problems as \emph{geodesically} convex problems over manifolds~\citep{zhang2016first,vishnoi2018geodesic}, 2) significant computational or performance gains in introducing manifolds constraints in some applications, e.g., training deep neural networks by enforcing hidden-layer weight matrices to belong to Stiefel manifolds \citep{arjovsky2016unitary,wisdom2016full,bansal2018can,xie2017all}, 3) new optimization problems that intrinsically involve manifold constraints~\citep{absil2019collection,hu2020brief,hosseini2020alternative,vandereycken2013low}.

Below, we present a policy evaluation problem in reinforcement learning which has two-level nested compositional form over Riemannian manifolds as \eqref{eq:main1}.

\vspace{-0.2cm}
\paragraph{Motivating example.} In reinforcement learning, finding the value function of a given policy is often referred to as policy evaluation problem \citep{dann2014policy,wang2017stochastic,sutton2018reinforcement}. Suppose there are $S$ states with the state variable $s\in\mR^d$ and denote the policy of interest by $\pi$. Let $V^\pi\in\mR^S$ denote the value function, where $V^\pi(s)$ represents the value of state $s$ under policy $\pi$. The Bellman equation that should be satisfied by the optimal policy is 
\begin{equation}
	V^\pi(s)=\mE [r_{s,s'}+\rho V^\pi(s')|s],\ \ \forall s, s' \in\{1,...,S\},
\end{equation}
where $\rho<1$ is a discount factor, $r_{s,s'}$ is the reward of transition from state $s$ to state $s'$, and the expectation is taken over all possible future states $s'$ conditioned on current state $s$ and the policy $\pi$. In the blackbox simulation environment, the transition and reward matrices are unknown, but can be sampled. Furthermore, a large number of states $S$ makes solving the Bellman equation directly impractical. Therefore, we model the value function with Gaussian basis functions (see \cite{bishop2006pattern}) and consider an iterative procedure to estimate the parameters. In particular, we assume that $V^\pi(s)\approx \sum_{i=1}^D w_i N(s;\mu_i,\Sigma_i)$, where $D\ll S$. We denote $\sum_{i=1}^D w_i N(s;\mu_i,\Sigma_i)$ by $\varphi(s;\Lambda)$, where $\Lambda\triangleq\{\{w_i\},\{\mu_i\},\{\Sigma_i\}\}_{i=1}^D$ is the set of all parameters.

The policy evaluation problem is formulated as
\begin{equation}\label{eq:policy_approx_problem}
	\min_{\Lambda} \sum_s [\varphi(s;\Lambda)-\sum_{s'}\mE[\hat{P}_{s,s'}](\mE[\hat{r}_{s,s'}]+\rho\varphi(s';\Lambda))]^2
\end{equation}
where the optimization variables are on the product of the Euclidean spaces and the manifolds of positive definite matrices. Problem \eqref{eq:policy_approx_problem} is an instance of \eqref{eq:main1} with $f(\cdot)=\norm{\cdot}^2$ and function $g$ defined as 
$
g: \mR^D\times \mR^{D\times d}\times\cM^D\longrightarrow \mR^{S}, 
\{\{w_i\},\{\mu_i\},\{\Sigma_i\}\}_{i=1}^D \longmapsto (g_1,...,g_{S})$,
\begin{equation}
g_s\triangleq\varphi(s;\Lambda)-\sum_{s'}\mE[\hat{P}_{s,s'}](\mE[\hat{r}_{s,s'}]+\rho\varphi(s';\Lambda)).
\end{equation}

\vspace{-0.2cm}
\subsection{Related work}
This work stems from two different lines of research: 1. Manifold optimization, 2. Stochastic compositional optimization in Euclidean setting which are briefly reviewed below.

\paragraph{Manifold optimization.} The stochastic gradient descent (SGD) method over manifolds was first studied in~\cite{bonnabel2013stochastic} which proved Riemannian SGD converges to a critical point of the problem. Under geodesic convexity, \cite{zhang2016first} developed the first global complexity of first-order methods (in general) and established $\cO(1/\epsilon^2)$ complexity to attain an $\epsilon$-optimal solution, i.e., $f(x^k)-f^*\leq \epsilon$, for Riemannian SGD. In the Euclidean setting, many variance reduction techniques have been proposed to improve the sample complexity of SGD~\citep{roux2012stochastic,johnson2013accelerating,xiao2014proximal,defazio2014saga,reddi2016stochastic,allen2016variance}. As a generalization of~\cite{johnson2013accelerating}, the Riemannian stochastic variance-reduced gradient descent (R-SVRG) method was developed in \cite{zhang2016riemannian}, establishing the linear rate for geodesically strongly convex functions. An extension of  Riemannian SVRG with computationally more efficient retraction and vector transport was developed in \cite{sato2019riemannian}. In~\cite{tripuraneni2018averaging}, authors adapted the Polyak-Ruppert iterate averaging over Riemannian manifolds~\citep{ruppert1988efficient,polyak1990new} reaching $\cO(1/k)$ rate. The Riemannian version of the stochastic recursive gradient method~\citep{nguyen2017stochastic} was proposed in~\cite{kasai2018riemannian}.~\cite{kasai2019riemannian} proposed the adaptive gradient method with the convergence rate of $\cO(\log(k)/\sqrt{k})$. 
 
 Besides the stochastic gradient methods, numerous deterministic algorithms for Euclidean unconstrained optimization~\citep{nocedal2006numerical,ruszczynski2011nonlinear} have also been generalized to Riemannian settings~\citep{absil2009optimization,udriste2013convex,boumal2020introduction}
- see, for instance, gradient descent~\citep{zhang2016first}, the Nesterov's accelerated method~\citep{liu2017accelerated,zhang2018towards,ahn2020nesterov,alimisis2020continuous}, proximal gradient method~\citep{huang2021riemannian}, Frank-Wolfe method~\citep{weber2019projection}, nonlinear conjugate gradient method~\citep{smith1994optimization,sato2015new}, BFGS~\citep{ring2012optimization}, Newton's method~\citep{adler2002newton,kasai2018riemannian}, trust-region method~\citep{absil2007trust} and cubic regularized method~\citep{zhang2018cubic,agarwal2018adaptive,zhang2020riemannian}.
 
 \paragraph{Stochastic compositional optimization.}
 The stochastic compositional optimization is closely related to the classical SGD and stochastic approximation (SA) methods~\citep{robbins1951stochastic,kiefer1952stochastic,kushner2003stochastic}. When the outer function $f$ is linear, the problem reduces to the standard stochastic non-compositional setting which has been extensively studied. While the problem is indeed compositional with two expected-value functions, \cite{ermoliev1976methods} applied a simple two-timescale SA scheme, showing its convergence under basic assumptions. 
 
The first nonasymptotic analysis for the stochastic compositional optimization appeared in~\cite{wang2017stochastic}, which uses two sequences of stepsizes in two different timescales: a slower one $\alpha_k$ for updating the optimization variable $\bx$ and a faster one $\beta_k$ for tracking the value of inner function. Their analysis requires $\alpha_k/\beta_k\to 0$ as $k\to\infty$. For problems with smooth and convex composition objective, their algorithm converge at a rate of $\cO(k^{-2/7})$, and with the rate of $\cO(k^{-4/5})$ in the strongly convex case, where $k$ represents the number of queries to the stochastic first-order oracles. The first finite-sample error bound is improved in~\cite{wang2016accelerating} to $\cO(k^{-4/9})$ for convex and nonconvex settings. While most methods rely on the two-timescale stepsizes, the single timescale algorithm was recently developed in~\cite{ghadimi2020single}, achieving the sample complexity of $\cO(1/\epsilon^2)$ to find an $\epsilon$-approximate stationary point. Furthermore, \cite{chen2021solving} proposed a single-loop loop algorithm, without any need for accuracy-dependent stepsize or increasing batch size, that can achieve the sample complexity of $\cO(1/\epsilon^2)$ as in \cite{chen2021solving} (or classic SGD for non-compositional problems). 
 
 In addition to the general stochastic optimization, the special setting with finite-sum structure recently gained popularity. The variants of the algorithm in \cite{wang2017stochastic} for finite-sum setting have been proposed in~\cite{lian2017finite,blanchet2017unbiased,devraj2019stochastic,lin2018improved,xu2021katyusha}. Furthermore, the stochastic compositional problem with certain nonsmooth component was investigated in~\cite{huo2018accelerated,zhang2019stochastic,zhang2019composite}. A key feature of these works is involving variance reduction techniques, which helps to achieve better performance for the finite-sum stochastic compositional problems. As these methods usually require increasing batch size, it is not possible to directly use them for general stochastic compositional problems.

 \subsection{Contributions}
The main contributions of this paper are as follows: 
\textbf{1)} The paper proposes algorithms to optimize the composition of two functions in the expectation form over Riemannian manifolds. To our knowledge, this is the first work opening the discussion for future works. The algorithms are motivated by the Riemannian extension of the ODE gradient flow presented in \cite{chen2021solving} for the unconstrained Euclidean setting. 
\textbf{2)} We provide the sample complexity of the proposed algorithms obtaining $\cO(1/\epsilon^2)$ to obtain $\epsilon$-approximate stationary solution, i.e., $\norm{\text{grad} f(x)}^2\leq \epsilon$, which is the same rate as Riemannian SGD for stochastic non-composition problems~\cite{hosseini2020alternative} or the algorithms in \cite{chen2021solving} and \cite{ghadimi2020single} for Euclidean composition problems.
\textbf{3)} We empirically verify the effectiveness of the proposed algorithm for two-level composition problems in the policy evaluation problem discussed in section~\ref{sec:intro}. 

\subsection{Preliminaries}
A Riemannnian manifold $(\cM,g)$ is a real smooth manifold $\cM$ equipped with a Riemannain metric $g$. The metric $g$ induces an inner product structure in each tangent space $T_{\bx} \cM$ associated with point $\bx\in\cM$. We denote the inner product of $\bu, \bv\in T_{\bx}\cM$ as $\fprod{\bu,\bv}_x$, and the norm of $\bu$ is defined as $\norm{\bu}=\sqrt{\fprod{\bu,\bu}_x}$. Given a smooth real-valued function $f$ on a Riemannian manifold $\cM$, Riemannian gradient of $f$ at $\bx$ is denoted by $\text{grad}f(\bx)$. We use $\nabla f$ to denote the gradient (or Jacobian) of a scalar (or vector) valued function $f$ in the Euclidean sense.

\begin{definition}
	Given manifolds $\cM$ and $\cM'$, the differential of a smooth map $F:\cM\to\cM'$ at $\bx$ is a linear operator $DF(\bx): \cT_{\bx}\cM\to\cT_{F(\bx)}\cM'$ defined by: 
	\begin{equation}
		DF(\bx)[v]=[t\mapsto F(c(t))], 
	\end{equation}
	where $c$ is a smooth curve on $\cM$ passing through $\bx$ at $t=0$ such that $v=[c]$. {Note that the bracket denotes the equivalence class of curve $c$, with the same ``velocity''-- see Definition~8.37 in \cite{boumal2020introduction} for more details.} In particular, when $\cM'$ is an embedded submanifold of a linear space $\cE$, identifying the tangent spaces of $\cM'$ to subspaces of $\cE$ and with $v=[c]$ a tangent vector at $\bx\in\cM$, we write 
	\begin{equation}\label{def:differential}
		DF(\bx)[v]=(F\circ c)'(0), 
	\end{equation} 
	where $F\circ c$ is seen as a map into $\cE$. 
\end{definition}

\begin{definition}[Riemannian gradient]\label{def:riemannian_gradient}
	Let $f: \cM\rightarrow\mR$ be a smooth function on a Riemannian manifold $\cM$. The Riemannian gradient of $f$ is the vector field $\text{grad}f$ on $\cM$ that $\forall \bx\in\cM, \forall v \in T_{\bx}\cM$ satisfies
	\begin{equation}
		Df(\bx)[v]=\fprod{v,\text{grad} f(\bx)}_{\bx},
	\end{equation}
	where $Df(\bx)$ is the differential of $f$ at $\bx\in\cM$.
\end{definition}

\begin{definition}[Adjoint of an operator]\label{def:adjoint}
	Let $\cE$ and $\cE'$ be two Euclidean spaces, with inner products $\fprod{\cdot,\cdot}_a$ and $\fprod{\cdot,\cdot}_b$ respectively. Let $\cA: \cE\mapsto\cE'$ be a linear operator. The adjoint of $\cA$ is a linear operator $\cA^*:\cE'\mapsto\cE$ and we have 
	\begin{equation}
		\forall \bu\in\cE, \bv\in\cE',\  \ \fprod{\cA(\bu),\bv}_b=\fprod{\bu,\cA^*(\bv)}_a.
	\end{equation}
\end{definition}

\section{Two-level Riemannian composition}\label{sec:algo}
We first characterize the Riemannian gradient of the composite function $F(\bx)$ 
\begin{align*}
	DF(\bx)[v]&= Df(g(\bx))[Dg(\bx)[v]]\\
	&=\fprod{\text{Proj}_{g(\bx)}\nabla f (g(\bx)),Dg(\bx)[v]}\\
	&=\fprod{(Dg(\bx))^*(\text{Proj}_{g(\bx)}\nabla f (g(\bx))),v}_{\bx}
\end{align*}
where $\text{Proj}_{g(\bx)}$ is the orthogonal projection onto $T_{g(\bx)}\cN$ follows from the embedding assumption of $\cN$ in the Euclidean space. 
According to Definition~\ref{def:riemannian_gradient}, the Riemannian gradient of $F(\bx)$, denoted by $\text{grad} F(\bx)$, is $(Dg(\bx))^*(\text{Proj}_{g(\bx)}\nabla f (g(\bx)))$. 

Following the standard SGD methodology applied on manifold optimization~\citep{robbins1951stochastic,bonnabel2013stochastic,zhang2016riemannian,zhang2016first,hosseini2020alternative}, a promising update is 
\begin{align}
	\bx^{k+1}=\text{Exp}_{\bx^k}(-\alpha(Dg_{\phi^k}(\bx^k))^*(\text{Proj}_{g(\bx^k)}\nabla f_{\xi^k} (g(\bx^k)))),
\end{align}
where $\phi^k$ and $\xi^k$ are samples drawn at iteration $k$.
Note that $Dg_{\phi^k}(\bx^k))^*(\text{Proj}_{g(\bx^k)}\nabla f_{\xi^k} (g(\bx^k))$ is an unbiased estimator of $\text{grad } F(\bx)$ but the exact evaluation of $g(\bx^k)$, i.e. $\mE_{\phi}[g_{\phi}(\bx^k)]$, is generally not attainable. Furthermore, the stochastic gradient is \emph{not} unbiased if one replaces $g(\bx^k)$ by its stochastic estimate $g_{\phi^k}(\bx^k)$. 
Therefore, the stochastic gradient method cannot be directly applied.

\subsection{Algorithm development motivated by ODE analysis}
Below, we provide the intuition behind our algorithmic design based on the ODE gradient flow which carefully extends the analysis in \cite{chen2021solving} to the Riemannian setting.

Let $t$ denote the time in this subsection. Consider the following ODE
\begin{equation}\label{algo_dev:ODE}
	\dot{\bx}(t)=-\alpha(Dg(\bx))^*\text{Proj}_{g(\bx)}\nabla f (\by(t)),
\end{equation}
for some $\alpha >0$. If we set $\by(t)=g(\btheta(t))$, then we have 
\begin{align*}
	\frac{d}{dt}f(g(\btheta(t)))=\fprod{(Dg(\bx))^*\text{Proj}_{g(\bx)}\nabla f (g(\btheta(t))),\dot{\btheta}(t)}_{\btheta(t)}
	=-\frac{1}{\alpha}\norm{\dot{\btheta}(t)}^2_{\btheta(t)}\leq 0.
\end{align*}
In this case, \eqref{algo_dev:ODE} describes a gradient flow that monotonically decreases $f(g(\btheta(t)))$.
However, we can not evaluate $g(\btheta(t))$ exactly. Instead, we can evaluate $\nabla f$ at $\by(t)\approx g(\btheta(t))$, and the introduced inexactness results in $f(g(\btheta(t)))$ loosing its monotonicity:

{\small
\begin{align*}
	\frac{d}{dt}f(g(\btheta(t)))
	&=-\frac{1}{\alpha}\norm{\dot{\btheta}(t)}^2_{\btheta(t)}+
	\fprod{(Dg(\bx))^*\text{Proj}_{g(\bx)}(\nabla f (g(\btheta(t)))-\nabla f(\by(t))),\dot{\btheta}(t)}_{\btheta(t)}\\
	&= -\frac{1}{\alpha}\norm{\dot{\btheta}(t)}^2_{\btheta(t)}+
	\fprod{\text{Proj}_{g(\bx)}(\nabla f (g(\btheta(t)))-\nabla f(\by(t))),Dg(\bx)[\dot{\btheta}(t)]}_{g(\btheta(t))}\\
	&\leq -\frac{1}{\alpha}\norm{\dot{\btheta}(t)}^2_{\btheta(t)}+ \norm{Dg(\bx)[\dot{\btheta}(t)]}_{g(\btheta(t))}
	\norm{\text{Proj}_{g(\bx)}(\nabla f (g(\btheta(t)))-\nabla f(\by(t)))}_{g(\btheta(t))} \\
	&\leq -\frac{1}{\alpha}\norm{\dot{\btheta}(t)}^2_{\btheta(t)}+L_f C_g\norm{g(\bx(t))-\by(t)}\cdot\norm{\dot{\bx}(t)}_{\bx(t)}\\
	&\leq -\frac{1}{2\alpha}\norm{\dot{\btheta}(t)}^2_{\btheta(t)}+\frac{\alpha C_g^2L_f^2}{2}\norm{g(\btheta(t))-\by(t))}^2.
\end{align*}
}
\normalsize
This motivates an energy function, 
\begin{equation}
	V(t)=f(g(\btheta(t)))+\norm{g(\btheta(t))-\by(t)}^2.
\end{equation}
We want $V(t)$ to monotonically decrease. By substitution, we have 
\small
\begin{align*}
	\dot{V}(t)
	&\leq -\frac{1}{2\alpha}\norm{\dot{\btheta}(t)}^2_{\btheta(t)}+\frac{\alpha C_g^2L_f^2}{2}\norm{g(\btheta(t))-\by(t))}^2+2\fprod{\by(t)-g(\btheta(t)),\dot{\by}(t)-Dg(\btheta(t))[\dot{\btheta}(t)]}\\
	&= -\frac{1}{2\alpha}\norm{\dot{\btheta}(t)}^2_{\btheta(t)}-(\beta-\frac{\alpha C_g^2L_f^2}{2})\norm{g(\btheta(t))-\by(t)}^2+2 \langle\by(t)-g(\btheta(t)), \dot{\by}(t)+\beta(\by(t)-g(\btheta(t))) \\
	&-Dg(\btheta(t))[\dot{\btheta}(t)]\rangle,
\end{align*}
\normalsize
where $\beta >0$ is a fixed constant. 
Following the maximum descent principle of $V(t)$, we are motivated to use the following dynamics
\begin{equation}\label{algo_dev:y}
	\dot{\by}(t)=-\beta(\by(t)-g(\btheta(t)))+Dg(\btheta(t))[\dot{\btheta}(t)]. 
\end{equation}
We approximate $Dg(\btheta(t))[\dot{\btheta}(t)]$ by the first-order Taylor expansion~\citep{boumal2020introduction}, i.e.,
\begin{equation}
	Dg(\btheta(t))[\dot{\btheta}(t)]\approx \gamma_k(g(\bx^k)-g(\bx^{k-1})),
\end{equation}
where $k$ is the discrete iteration index, and $\gamma_k$ is the weight controlling the approximation. 

With the insights gained from~\eqref{algo_dev:ODE} and~\eqref{algo_dev:y}, we propose the following stochastic update, which serve as the main components in Algorithm~\ref{alg:SCSC1}. 
\begin{align}
	\by^{k+1}=&\by^k-\beta_k(\by^k-g_{\phi^k}(\bx^k)) +\gamma_k(g_{\phi^k}(\bx^k)-g_{\phi^k}(\bx^{k-1})),
\end{align}
\begin{align}
	\bx^{k+1}= \text{Exp}_{\bx^k}(-\alpha (Dg_{\phi^k}(\bx^k))^*\text{Proj}_{g_{\phi^k}(\bx^k)}\nabla f_{\xi^k}(\by^{k+1})).
\end{align}

\begin{algorithm}[H]
\caption{R-SCGD for two-level problem}
\label{alg:SCSC1}
\begin{algorithmic}[1]
\REQUIRE $\bx^0$, $\by^0$, constant sequences $\{\alpha_k\}$, $\{\beta_k\}$, $\{\gamma_k\}$
\FOR{$k=0,...,K-1$} 
\STATE Randomly sample $\phi^k$ and $\xi^k$  
\STATE Update inner function estimate $\by^{k+1}$ using
$$\by^{k+1}=\by^k-\beta_k(\by^k-g_{\phi^k}(\bx^k))+\gamma_k(g_{\phi^k}(\bx^k)-g_{\phi^k}(\bx^{k-1}))$$

\STATE Update Riemannian gradient of the composition function as
\vspace{-0.2cm}
\begin{align*}
\eta^{k+1}=(Dg_{\phi^k}(\bx^k))^* \text{Proj}_{g_{\phi^k}(\bx^k)}\nabla f_{\xi^k}(\by^{k+1})
\end{align*}
\STATE Update $\bx^{k+1}=\text{Exp}_{\bx^k}(-\alpha_k \eta^{k+1})$
\ENDFOR

\end{algorithmic}
\end{algorithm}

\subsection{Iteration complexity of the two-level R-SCGD method}\label{analysis_for_bilevel}
With the insights gained from the continuous-time Lyapunov function, our analysis in this section essentially builds on the following discrete-time Lyapunov function
\begin{equation}\label{lyapunov}
	\cV^k:=F(\bx^k)-F(\bx^*)+\norm{g(\bx^{k-1})-\by^k}^2,
\end{equation}
where $\bx^*$ is (one of) the optimal solution(s) of the problem~\eqref{eq:main1}.

\begin{assumption}\label{assump:smooth_obj}
	Function $F(\bx):\cM\to\mR$ is geodesically $L_g$-smooth, i.e., $\forall (\bx,\bs)\in\cT$,
	\begin{equation}
	  \norm{P^{-1}_s\text{grad }F(\text{Exp}_{\bx}(s))-\text{grad }F(\bx)}\leq L_g\norm{s},
	\end{equation}
	where $\cT\subseteq T\cM$ is the domain of \text{Exp} and $P_s$ denotes parallel transport along $\gamma(t)=\text{Exp}_{\bx}(ts)$ from $t=0$ to $t=1$.
\end{assumption}

\begin{assumption}\label{assumption:smooth_f}	
	Function $f_{\xi}: \cE\to\mR$ is $L_f$-smooth, i.e., for all $\by,\by'\in\cE$,
	\begin{equation}\label{f_lip}
		\norm{\nabla f_{\xi}(\by)-\nabla f_{\xi}(\by')}\leq L_f \norm{\by-\by'}.
	\end{equation}
\end{assumption}
\begin{assumption}\label{assump:unbiased}
	Random sample oracle of function value $g_{\phi}(\bx)$ is an unbiased estimator of $g(\bx)$ and has bounded variance, i.e., 
	\begin{align}
	\begin{split}
		\mE[g_{\phi}(\bx)]=g(\bx), \\
		\mE[\norm{g_{\phi}(\bx)-g(\bx)}^2]\leq V_g^2.
	\end{split}
	\end{align}
\end{assumption}

\begin{assumption}\label{assump:oracle}
	The chain rule holds in expectation, i.e. 
	\begin{align*}
		\mE[(Dg_{\phi}(\bx))^*\text{Proj}_{g_{\phi}(\bx)}\nabla f_{\xi}(g(\bx))]=(Dg(\bx))^*\text{Proj}_{g(\bx)}\nabla f(g(\bx)).	
	\end{align*}
\end{assumption}

\begin{remark}\label{remark:1}
	If the random sample oracle of derivatives (gradients) satisfies
	\begin{align}
	\begin{split}
		\mE[Dg_{\phi}(\bx)[\eta]]=Dg(\bx)[\eta], \\
		\mE[\nabla f_{\xi}(\bx)]=\nabla f(\bx), 
	\end{split}
	\end{align}
	and with the independence between $\phi$ and $\xi$, Assumption~\ref{assump:oracle} holds. See lemma below. 
\end{remark}

\begin{lemma}\label{lemma:chainrule}
	If the random sample oracle for derivatives (gradients) satisfies
	\begin{equation}
		\mE[Dg_{\phi}(\bx)[\eta]]=Dg(\bx)[\eta], \text{and } \mE[\nabla f_{\xi}(\bx)]=\nabla f(\bx), 
	\end{equation}
	and random variable $\phi$ is independent of $\xi$, then
	\begin{equation}
		\mE[(Dg_{\phi}(\bx))^*\text{Proj}_{g_{\phi}(\bx)}\nabla f_{\xi}(g(\bx))]=(Dg(\bx))^*\text{Proj}_{g(\bx)}\nabla f(g(\bx)).	
	\end{equation}
\end{lemma}

\begin{assumption}\label{assumption:bounded1}
	The stochastic gradients of $f_{\xi}$ and $g_{\phi}$ are bounded in expectation, i.e.,
	\begin{align}
	\begin{split}
		\mE [\norm{Dg_{\phi}(\bx)[\cdot]}^2_{op}]\leq C_g^2, \\
		\mE [\norm{\nabla f_{\xi}(\bx)}^2]\leq C_f^2.
	\end{split}
	\end{align}
\end{assumption}

The Assumptions~\ref{assump:smooth_obj} and~\ref{assumption:smooth_f} require that the objective function and the outer function have lipshitz continuous gradients. The Assumptions~\ref{assump:unbiased}-\ref{assumption:bounded1} require the stochastic oracles satisfy certain unbiasedness and second-moment boundedness, which are typical assumptions for stochastic methods. When the manifold $\cM$ in problem~\eqref{eq:main1} degenerates to a linear space, then the presented assumptions are equivalent to the assumptions in stochastic compositional optimization in the Euclidean setting (\cite{wang2017stochastic,wang2016accelerating,zhang2019stochastic,chen2021solving}). 

\begin{lemma}\label{lemma:variance_tracking}
	Consider $\cF^k$ as the collection of random variables, i.e., $$\cF^k:=\{\phi^0,...,\phi^{k-1},\xi^0,...,\xi^{k-1}\}.$$ Suppose that the Assumptions~\ref{assump:unbiased} and~\ref{assumption:bounded1} hold, and $\by^{k+1}$ is generated by running R-SCGD in Algorithm~\ref{alg:SCSC1} conditioned on $\cF^k$, then the mean square error of $\by^{k+1}$ satisfies 
\begin{align}
	\mE[\norm{\by^{k+1}-g(\bx^k)}^2|\cF^k] \leq (1-\beta_k)\norm{\by^k-g(\bx^{k-1})}^2 +2\beta_k^2V_g^2+B_0\alpha_k^2C_g^4C_f^2,
\end{align}
	where the constant is defined as $B_0:=2(1-\beta_k)^2+2\gamma_k^2+\frac{(1-\beta_k-\gamma_k)^2}{\beta_k}$.
\end{lemma}

\begin{theorem}[Two-level R-SCGD]\label{mainthm1}
Under the Assumption~\ref{assump:smooth_obj} to~\ref{assumption:bounded1}, if we choose the stepsizes as $\alpha_k=\frac{2\beta_k}{C_g^2L_f^2}=\frac{1}{\sqrt{K}}$ and 
 $\gamma_k=1-t_k\beta_k$, such that $t:=\sup\{|t_k|\}$ is finite, the iterates $\{\bx^k\}$ of SCSC in Algorithm~\ref{alg:SCSC1} satisfy
	\begin{equation}
		\frac{\sum_{k=0}^{K-1}\mE[\norm{\text{grad }F(\bx^k)}^2]}{K}\leq \frac{2\cV^0+2B_1}{\sqrt{K}}
	\end{equation}
where the constant is defined as $B_1:=\frac{L_F}{2}C_g^2C_f^2+C_g^4L_f^4V_g^2+[2+6(t+1)^2]C_g^4C_f^2$. 
\end{theorem}

\section{Multi-level Riemannian composition}\label{sec:multi-level}
In this section, we consider the following multi-level compositional problem
\begin{align}\label{eq:main2}
\min_{\bx\in \cM} F(\bx) &= f_N(f_{N-1}(...f_1(\bx)...)),
\text{ where } \quad f_n(\cdot)=\mE_{\theta_n}[f_n(\cdot\ ;\theta_n)]
\end{align}
and $\bx\in\cM$ is the optimization variable lying on a smooth Riemannian manifold, function $f_1:\cM\to\cN\subseteq\mR^{d_2}$, $f_n:\mR^{d_n}\to\mR^{d_{n+1}}$, $n=2,3,...,N$ with $d_{N+1}=1$, are smooth but possibly nonconvex functions, and $\theta_1,...,\theta_N$ are independent random variables.

As a generalization to Algorithm~\ref{alg:SCSC1}, the Algorithm~\ref{alg:SCSC2} is proposed to tackle the problem~\eqref{eq:main2}, where we use $\by_1$,..., $\by_{N-1}$ to track the function values $f_1$,..., $f_{N-1}$. 

\begin{algorithm}[H]
\caption{R-SCGD for multi-level problem}
\label{alg:SCSC2}
\begin{algorithmic}[1]
\REQUIRE $\bx^0$, $\by^0$, constant sequences $\{\alpha_k\}$, $\{\beta_k\}$, $\{\gamma_k\}$
\FOR{$k=0,...,K-1$}
\STATE Randomly sample $\theta^k_1$
\STATE $\by_1^{k+1}=\by_1^k-\beta_k(\by_1^k-f_1(\bx^k;\theta_1^k))+\gamma_k(f_1(\bx^k;\theta_1^k)-f_1(\bx^{k-1};\theta_1^k))$ 
\FOR{$n=2,...,N-1$}
\STATE Randomly sample $\theta^k_n$
\STATE $\by^{k+1}_{n}=\by_n^k-\beta_k(\by^k_n-f_n(\by^{k+1}_{n-1};\theta^k_n))+\gamma_k(f_n(\by^{k+1}_{n-1};\theta^k_n)-f_n(\by^{k}_{n-1};\theta^k_n))$
\ENDFOR
\STATE Randomly sample $\theta^k_N$, and update Riemannian gradient of the composition as
\begin{align*}
\eta^{k+1}=(Df_1(\bx^k;\theta^k_1))^*
\text{Proj}_{f_1(\bx^k;\theta^k_1)}\nabla f_{2}(\by^{k+1}_1;\theta^k_2) \cdots
\nabla f_{N}(\by^{k+1}_{N-1};\theta^k_N)
\end{align*}
\STATE Update $\bx^{k+1}=\text{Exp}_{\bx^k}(-\alpha_k \eta^{k+1})$
\ENDFOR
\end{algorithmic}
\end{algorithm}

\subsection{Iteration complexity of the multi-level R-SCGD method}
Similar to the analysis in Section~\ref{analysis_for_bilevel} for the two-level R-SCGD, we will first quantify the error between the estimated function value $\by_n$ and exact function value $f_n(\by_{n-1})$, and then establish the complexity rate based on the generalized Lyapunov function
\begin{equation}\label{lyapunov:multi}
	\cV^k:=F(\bx^k)-F(\bx^*)+\sum_{n=1}^{N-1}\norm{\by_n^k-f_n(\by_{n-1}^k)}^2
\end{equation}
where $\bx^*$ is (one of) the optimal solution(s) of the problem~\eqref{eq:main2} and we will use $\by_0^k$ and $\bx^k$ interchangeably. Specifically, we enforce $\gamma^k=1-\beta^k$, which makes the presentation simpler.

The smoothness of the objective function is still needed. To analyze the multi-level R-SCGD algorithm, Assumptions~\ref{assumption:smooth_f}-\ref{assumption:bounded1} are revised as discussed below.

The smoothness of the objective function (Assumption~\ref{assump:smooth_obj}) is needed to establish Theorem~\ref{mainthm2}. Assumptions~\ref{assumption:smooth_f} to~\ref{assumption:bounded1} are revised for the multi-level setting as follows.

\begin{assumption}\label{assumption:smooth_f_multi}	
	For any $n$, function $f_{n}: \mR^{d_n}\to\mR^{d_{n+1}}$ is $L_n$-smooth, i.e., for all $\by,\by'\in\mR^{d_n}$,
	\begin{equation}\label{f_lip}
		\norm{\nabla f(\by;\theta_n)-\nabla f(\by';\theta_n)}\leq L_n \norm{\by-\by'}.
	\end{equation}
\end{assumption}

\begin{assumption}\label{assump:unbiased_multi}
	Random sample oracle of function value $\{f_n\}$ is an unbiased estimator and has bounded variance, i.e., 
	\begin{equation}
		\mE[f_n(\by_{n-1};\theta_n)]=f_n(\by_{n-1}), \text{and } \mE[\norm{f_n(\by_{n-1};\theta_n)-f_n(\by_{n-1})}^2]\leq V_n^2,\ \ n=1,...,N.
	\end{equation}
\end{assumption}

\begin{assumption}\label{assump:oracle_multi}
	The chain rule holds in expectation, i.e.,
	\begin{align}
		&\mE[(Df_1(\bx;\theta_1))^*\text{Proj}_{f_1(\bx;\theta_1)}\nabla f_2(f^{(1)}(\bx);\theta_2)\cdots\nabla f_N(f^{(N-1)}(\bx);\theta_N)]\\
		&=(Df_1(\bx))^*\text{Proj}_{f_1(\bx)}\nabla f_2(f^{(1)}(\bx))\cdots\nabla f_N(f^{(N-1)}(\bx)),
	\end{align}
	where $f^{(n)}(\bx)$ is defined as $f^{(n)}(\bx):=f_n\circ f_{n-1}\circ \cdots \circ f_1(\bx)$.
\end{assumption}

\begin{remark}
	Random sample oracles of derivatives (gradients) satisfy
	\begin{equation}
		\mE[Df_1(\bx;\theta_1)[\eta]]=Df_1(\bx)[\eta], \text{and } \mE[\nabla f_n(\by_{n-1};\theta_n)]=\nabla f_n(\by_{n-1}), \ \ n=2,...,N.
	\end{equation}
	and with the independence between $\{\theta_n,n=1,...,N\}$, the Assumption~\ref{assump:oracle_multi} holds. We omit the proof since it directly follows the proof in Lemma~\ref{lemma:chainrule}.
\end{remark}

\begin{assumption}\label{assumption:bounded1_multi}
	The stochastic gradients of $\{f_{n}\}$ are bounded in expectation, i.e.,
	\begin{equation}
		\mE [\norm{Df_1(\bx;\theta_1)[\cdot]}^2_{op}]\leq C_1^2, and\ \ \mE [\norm{\nabla f_{n}(\by_{n-1};\theta_n)}^2]\leq C_n^2, n=2,...,N.
	\end{equation}
\end{assumption}

\begin{lemma}\label{lemma:variance_tracking_multi}
	Consider $\cF^k$ as the collection of random variables, i.e. $\{\theta_i^j:1\leq i\leq N, 0\leq j\leq k-1\}$. Suppose that the Assumptions~\ref{assump:unbiased_multi} and ~\ref{assumption:bounded1_multi} hold, and $\by^{k+1}$ is the sequence generated by Algorithm~\ref{alg:SCSC2}, then the mean square error of $\by_n^{k+1}$ satisfies 
	\begin{equation}\label{result_lemma_multi}
	\begin{split}
		\mE[\norm{\by^{k+1}_n-f_n(\by^{k+1}_{n-1})}^2|\cF^k] & \leq
		 (1-\beta_k)\norm{\by^k_n-f_n(\by^k_{n-1})}^2+[2(1-\beta_k)^2+\frac{(1-\beta_k-\gamma_k)^2}{\beta_k} \\
		 &+2\gamma_k^2]C_n^2\norm{\by^k_{n-1}-\by^{k+1}_{n-1}}^2+2\beta_k^2V_n^2.
	\end{split}
	\end{equation}
\end{lemma}

\begin{theorem}[Multi-level R-SCGD]\label{mainthm2}
	Under the Assumption~\ref{assump:smooth_obj} and 6-9 (from the appendix), if we choose the stepsize $\alpha_k=\frac{2\beta_k}{\sum_{n=1}^{N-1}A_n^2}=\frac{1}{\sqrt{K}}$, the iterates $\bx^k$ of the multi-level R-SCGD in Algorithm~\ref{alg:SCSC2} satisfy
	\begin{equation}
		\frac{\sum_{k=0}^{K-1}\mE[\norm{\text{grad }F(\bx^k)}^2]}{K}=\cO(\frac{1}{\sqrt{K}}),
	\end{equation}
	where 
	$A_n$ and the constant in $\cO(\frac{1}{\sqrt{K}})$ are specified in the proof.
\end{theorem}

\vspace{-0.2cm}
\section{Numerical studies}
This section conducts numerical experiments over problem~\eqref{eq:policy_approx_problem} to compare the proposed algorithm with the straightforward implementation of the Riemannian SGD method. The code is written in MATLAB and uses the MANOPT package~\cite{boumal2014manopt}. All the studies are run on a laptop with a 1.4 GHz Quad-Core Intel Core i5 CPU and 8 GB memory.

A $9\times 9$ grid over the state space, i.e., $d=2$, is generated first. Next, we fix the number of basis functions to five, randomly generate the true parameters $\{w_i\}$, $\{\mu_i\}$, $\{\Sigma_i\}$ and generate the true value function, represented by a column vector $\bv$. We assume $\Sigma_i$ is unique across the basis functions and we fix $\{w_i\}$ and $\{\mu_i\}$ to their true value and optimize~\eqref{eq:policy_approx_problem} over $\Sigma$ belonging to the symmetric positive definite manifold. Next, we find the ``true'' transition matrix $P$ and reward matrix $r$ based on the true parameters such that the Bellman equation holds. The true transition and reward matrices are added with zero mean Gaussian noises, and the transition matrix is further normalized to be doubly stochastic. These two matrices are used as the outputs of the stochastic oracle at each iteration of the algorithm. We also generate an $11\times 11$ grid over the state space and follow a similar simulation procedure. 

\begin{figure}[H]
    \centering
    \begin{subfigure}[b]{0.3\textwidth}
        \centering
        \includegraphics[width=\textwidth]{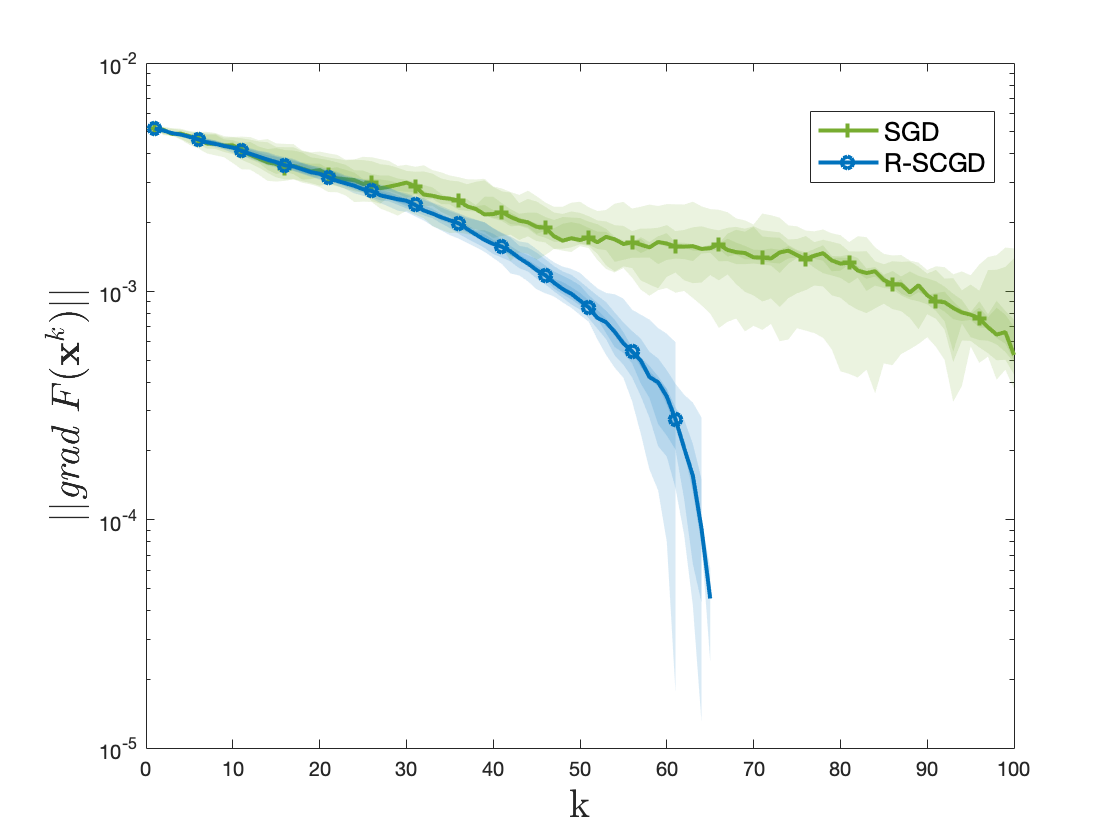}
    \end{subfigure}
    ~ 
    \begin{subfigure}[b]{0.3\textwidth}
        \centering
        \includegraphics[width=\textwidth]{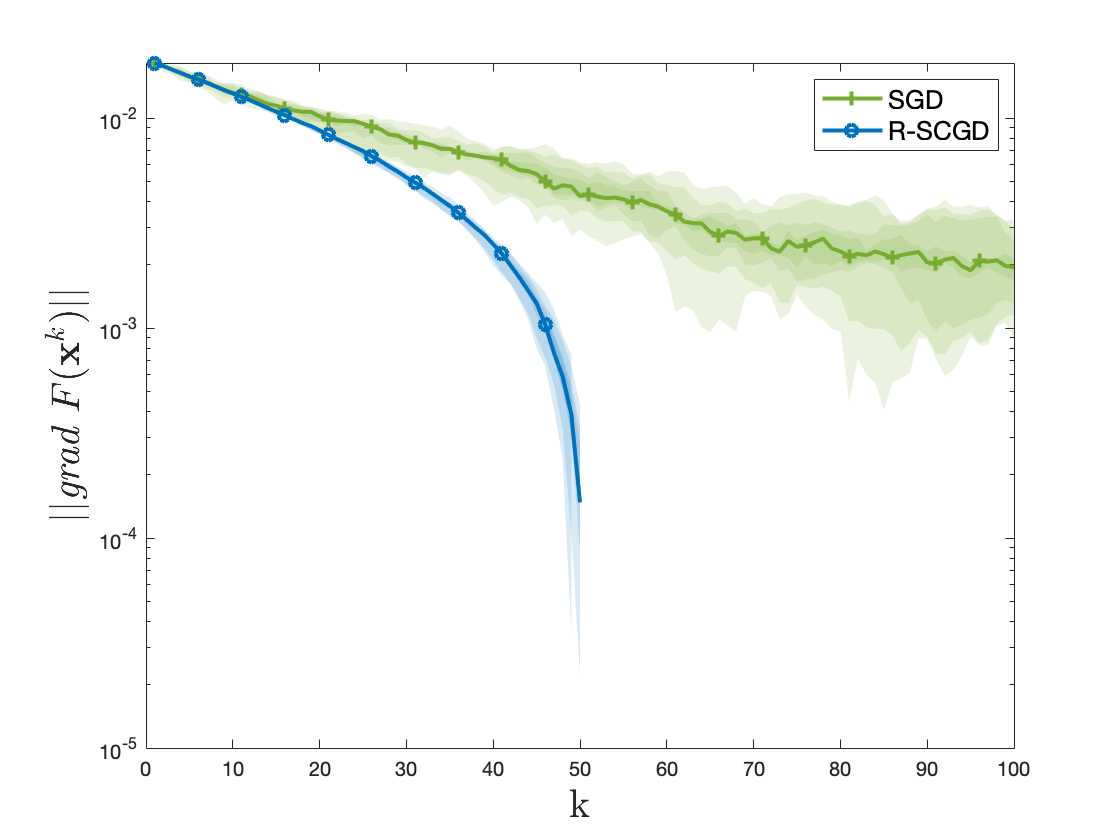}
    \end{subfigure}

    \begin{subfigure}[b]{0.3\textwidth}
        \centering
        \includegraphics[width=\textwidth]{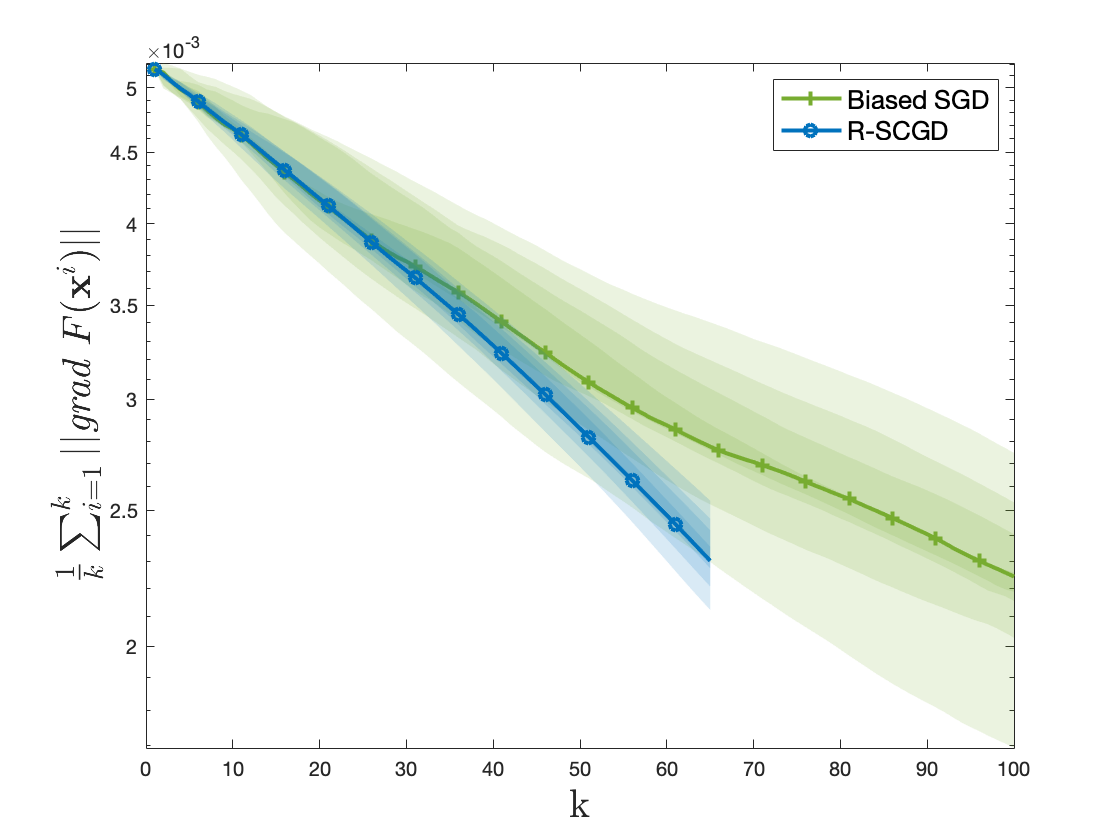}
    \end{subfigure}
    ~ 
    \begin{subfigure}[b]{0.3\textwidth}
        \centering
        \includegraphics[width=\textwidth]{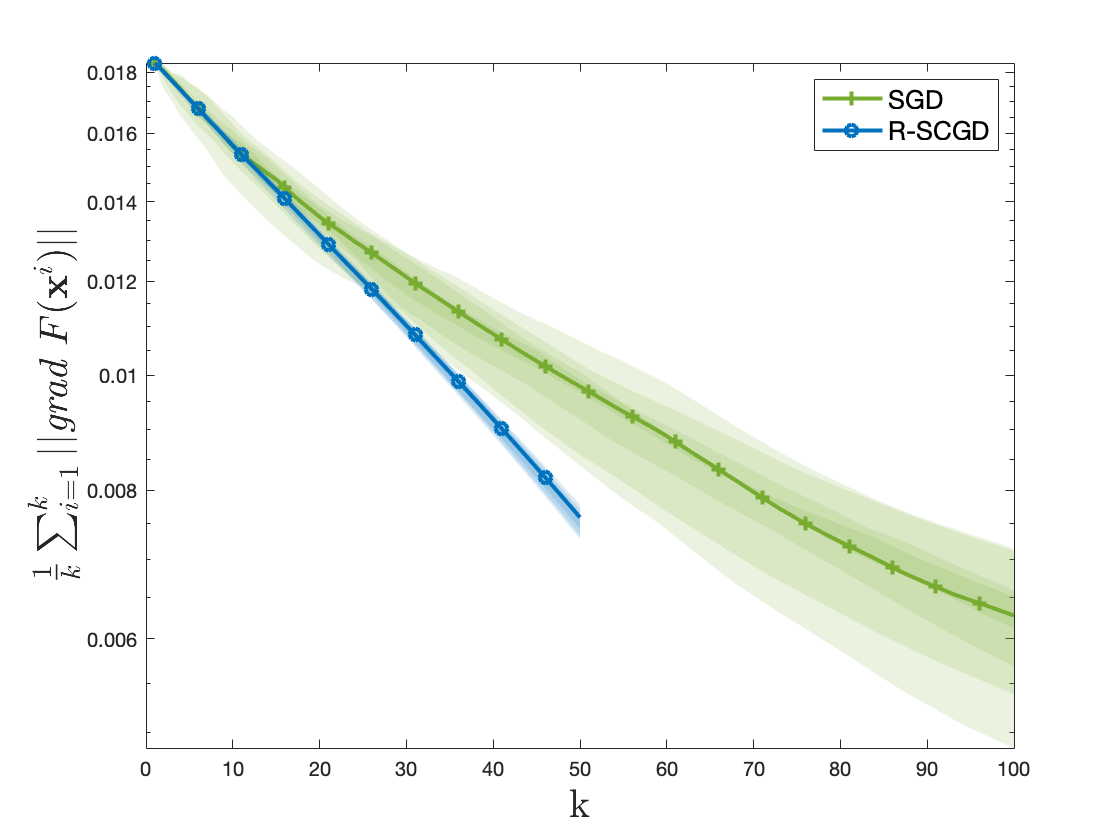}
    \end{subfigure}
    
    \caption{Performance of the proposed R-SCGD algorithm compared to the Riemannian SGD, for the norm of the Riemannian gradient. (\textbf{Left}) The first setting with 81 states. (\textbf{Right}) The second setting with 121 states.}
    \label{fig:grad_and_grad_avg}
\end{figure}

\begin{figure}[H]
    \centering
    \begin{subfigure}[b]{0.3\textwidth}
        \centering
        \includegraphics[width=\textwidth]{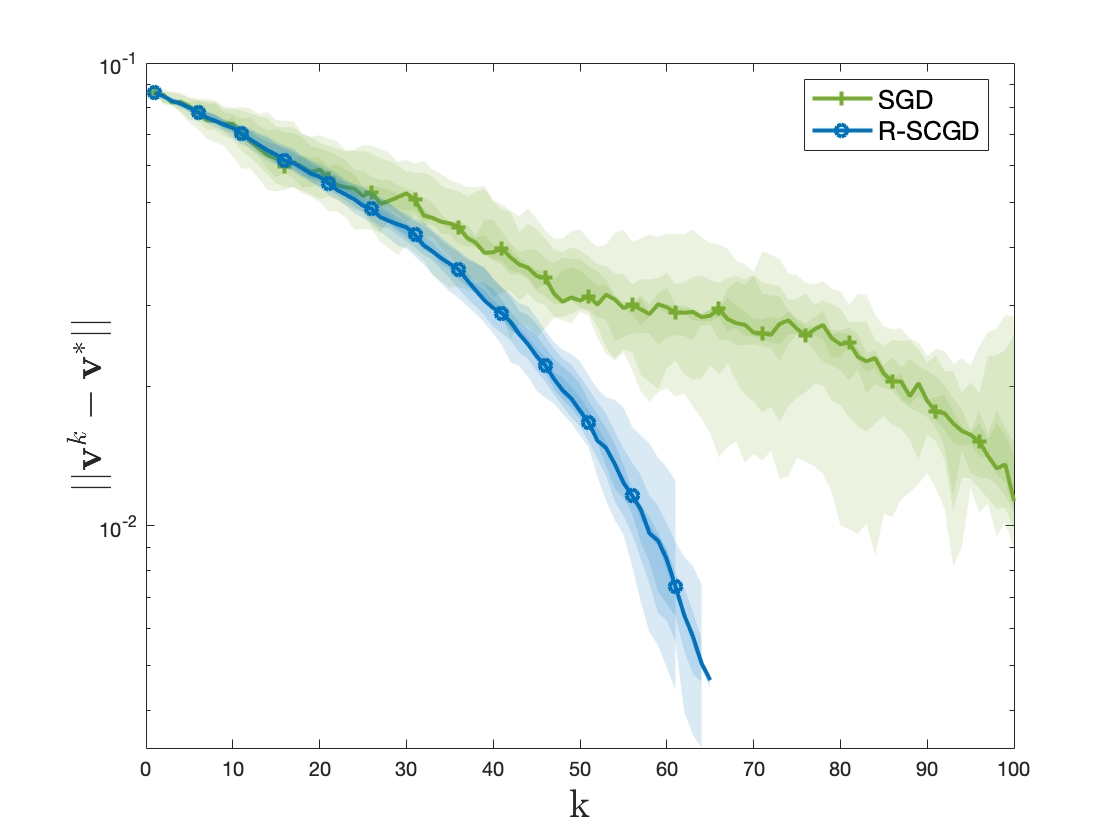}
    \end{subfigure}
    ~ 
    \begin{subfigure}[b]{0.3\textwidth}
        \centering
        \includegraphics[width=\textwidth]{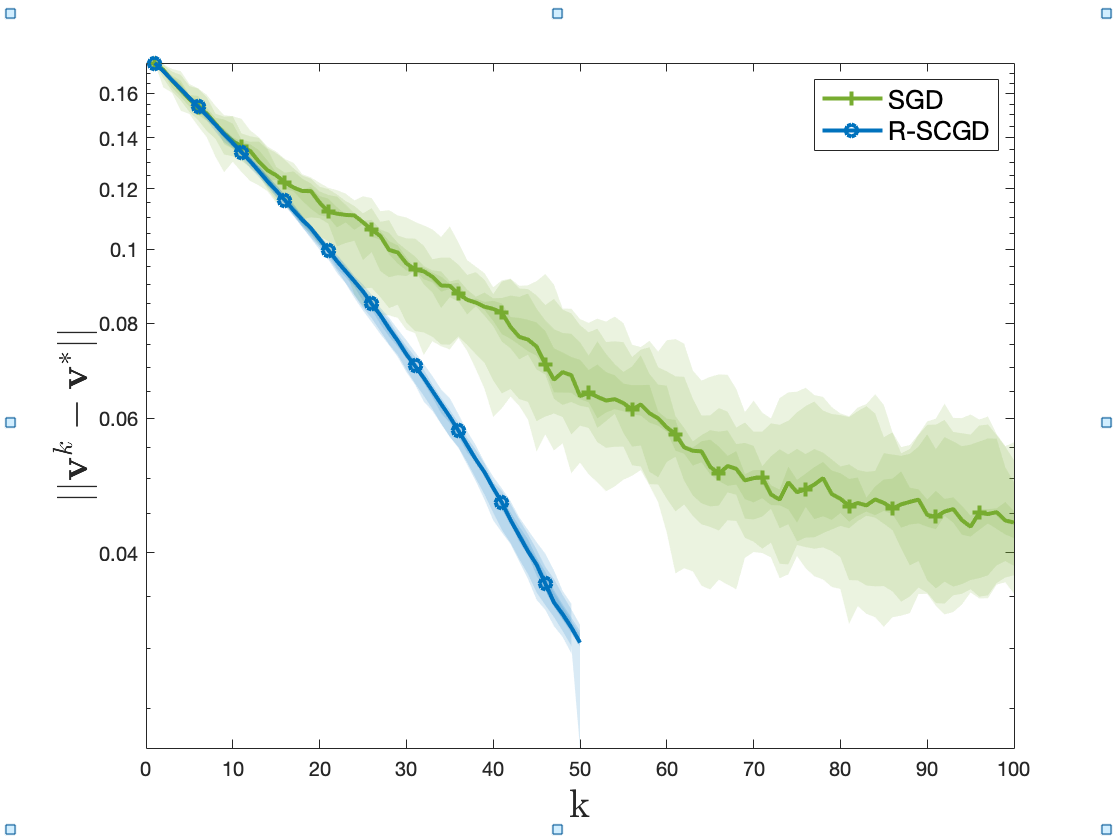}
    \end{subfigure}

    \begin{subfigure}[b]{0.3\textwidth}
        \centering
        \includegraphics[width=\textwidth]{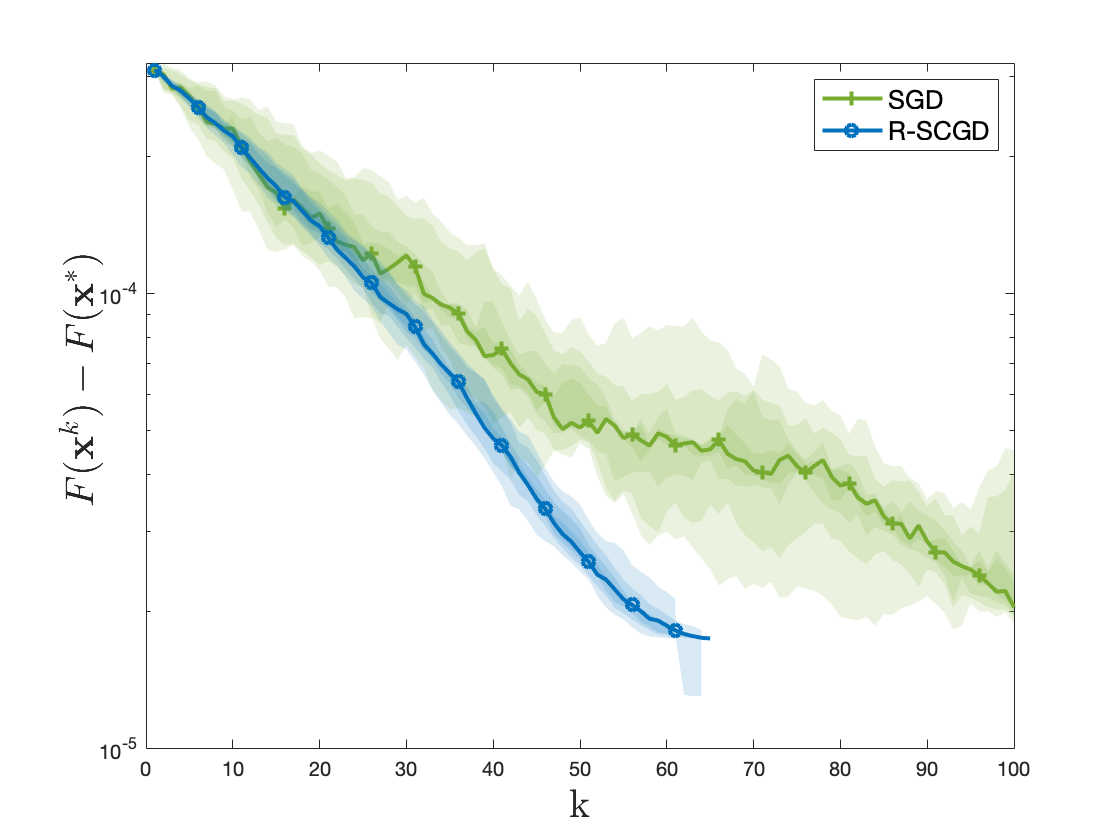}
    \end{subfigure}
    ~ 
    \begin{subfigure}[b]{0.3\textwidth}
        \centering
        \includegraphics[width=\textwidth]{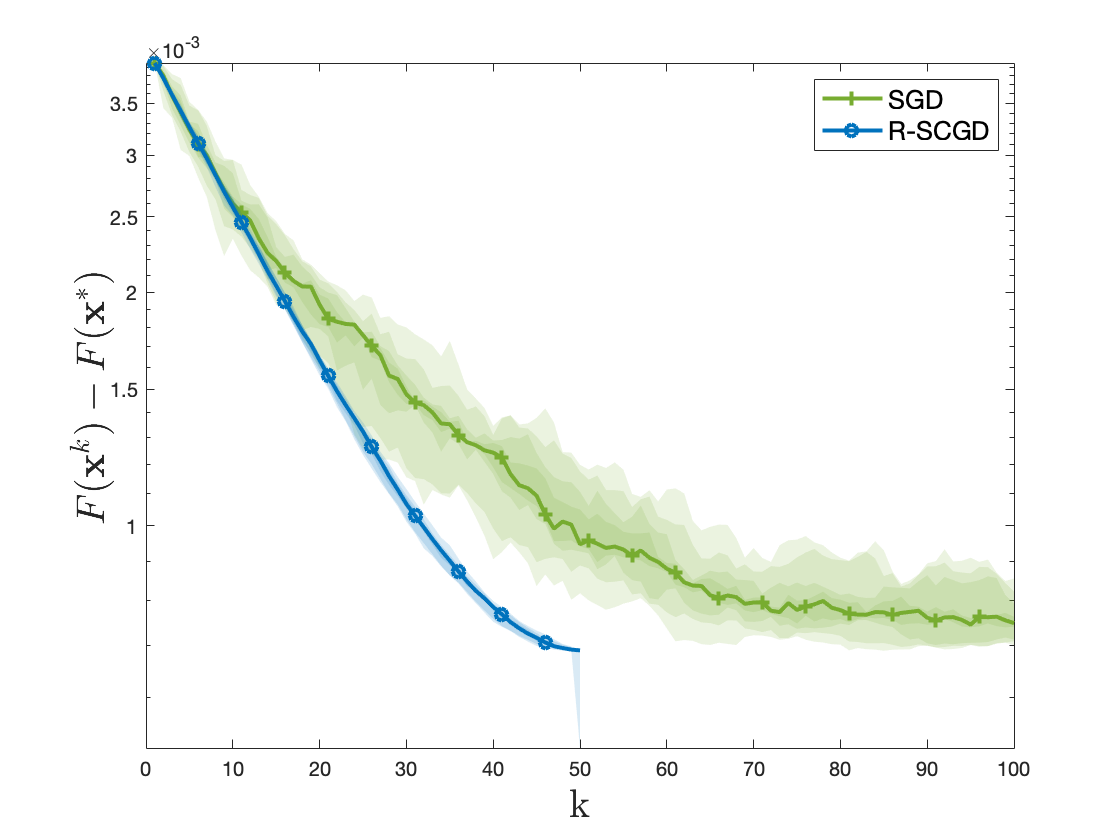}
    \end{subfigure} 
    \caption{Performance of the proposed R-SCGD algorithm compared to the Riemannian SGD. Left and right plots illustrate the 81- and 121-state settings, respectively. \textbf{(Top)} Top plots show the inner function approximation bias \textbf{(Bottom)} Bottom plots show the function value gap.
    }
    \label{fig:v_and_F}
\end{figure}

Figure~\ref{fig:grad_and_grad_avg} and~\ref{fig:v_and_F} compare the proposed method with the Riemannian SGD~\cite{bonnabel2013stochastic,zhang2016first} over 10 replicates for each scenario. The shades provide the percentile information based on the replicates. More specifically, Figure~\ref{fig:grad_and_grad_avg} presents the decreasing norm of the Riemannian gradients. The upper two plots select the last iterates as the output while the bottom two plots show the ergodic average of the iterates, on which the analysis is based. On those plots, the proposed R-SCGD method shows at least a linear convergence rate.

Though the test problem is not geodesically convex~\cite{boumal2020introduction}, the decreasing approximate bias and objective function value are illustrated in Figure~\ref{fig:v_and_F}. The results show better performance of the proposed R-SCGD compared to the biased Riemannian SGD for the composition problem.

\section{Conclusion}
We present the Riemannian stochastic compositional gradient method (R-SCGD) to solve the composition of two or multiple functions involving expectations over Riemannian manifolds. 
The proposed algorithm which is motivated by the Riemannian gradient flow approximates the inner function(s) value(s) using a moving average corrected by first-order information and its parameter(s) is in the same timescale as the stepsize of the variable update.
We established the sample complexity of $\cO(1/\epsilon^2)$ for the proposed algorithms to obtain $\epsilon$-approximate stationary solution, i.e., $\norm{\text{grad} f(x)}^2\leq \epsilon$.
We empirically studied the performance of the proposed algorithm over a policy approximation example in reinforcement learning.

\bibliography{refs.bib}


\vspace{0.2cm}
\section*{Appendices}
We provide the proofs related to the complexity analysis of the two-level and multi-level R-SCGD algorithms. 
The assumptions for the analysis of the multi-level R-SCGD algorithm, which are omitted from the body of the paper, are stated in Section~\ref{sec:supp_multi-level}.

\addcontentsline{toc}{section}{Appendices}
\renewcommand{\thesubsection}{\Alph{subsection}}

\vspace{-0.2cm}
\subsection{Complexity analysis of the two-level R-SCGD} \label{sec:supp_two-level}

\subsubsection{Proof of Lemma~\ref{lemma:chainrule}}
\begin{proof}
	It suffices to show that 
	\begin{equation}
		\fprod{\mE[(Dg_{\phi}(\bx))^*\text{Proj}_{g_{\phi}(\bx)}\nabla f_{\xi}(g(\bx))],b_i}_{\bx}=\fprod{(Dg(\bx))^*\text{Proj}_{g(\bx)}\nabla f(g(\bx)),b_i}_{\bx},
	\end{equation}
	where $\fprod{\cdot,\cdot}_{\bx}$ is the inner product defined on the tangent plane of the manifold $\cM$ at $\bx$, i.e. $T_{\bx}\cM$, and $\{b_i\}$ is the basis of $T_{\bx}\cM$. We have 
	\begin{align}
		&\fprod{\mE[(Dg_{\phi}(\bx))^*\text{Proj}_{g_{\phi}(\bx)}\nabla f_{\xi}(g(\bx))],b_i}_{\bx}\\
		&\labelrel={myeq:1}\mE\fprod{\text{Proj}_{g_{\phi}(\bx)}\nabla f_{\xi}(g(\bx)),Dg_{\phi}(\bx)[b_i]}_{g_{\phi}(\bx)}\label{2}\\
		&\labelrel={myeq:2}\mE\fprod{\nabla f_{\xi}(g(\bx)),Dg_{\phi}(\bx)[b_i]}\label{3}\\
		&\labelrel={myeq:3}\fprod{\nabla f(g(\bx)),Dg(\bx)[b_i]}\\
		&\labelrel={myeq:4}\fprod{(Dg(\bx))^*\text{Proj}_{g(\bx)}\nabla f(g(\bx)),b_i}_{\bx}.
	\end{align}
	The equality~\eqref{myeq:1} comes from the linearity of the inner product and the definition of adjoint operator (see Definition~\ref{def:adjoint}). For the equality~\eqref{myeq:2}, since we assume the manifold $\cN$ is embedded in the Euclidean space $\cE$, the inner product is the induced metric inherited from the Euclidean metric~\citep{boumal2020introduction} and the orthogonal projection operator can be erased. The equality~\eqref{myeq:3} follows from the independence between random variables $\xi$ and $\phi$. Finally, the last equality~\eqref{myeq:4} follows from the definition of adjoint operator.
\end{proof}

\subsubsection{Proof of Lemma~\ref{lemma:variance_tracking}}
\begin{proof}
Under the update rule
\begin{equation*}
\by^{k+1}=\by^k-\beta_k(\by^k-g_{\phi^k}(\bx^k))+\gamma_k(g_{\phi^k}(\bx^k)-g_{\phi^k}(\bx^{k-1})),
\end{equation*}
we have 
\begin{align}\label{lemmaproof:1}
	\by^{k+1}-g(\bx^k)=(1-\beta_k)(\by^k-g(\bx^{k-1}))+(1-\beta_k)T_{1,k}+\beta_kT_{2,k}+\gamma_kT_{3,k},
\end{align}
where $T_{1,k}:=g(\bx^{k-1})-g(\bx^k)$, $T_{2,k}:=g_{\phi^k}(\bx^{k})-g(\bx^k)$, $T_{3,k}:=g_{\phi^k}(\bx^{k})-g_{\phi^k}(\bx^{k-1})$.
Conditioned on $\cF^k$, taking expectation over $\phi^k$, we have 
\begin{align}
	\mE[(1-\beta_k)T_{1,k}+\beta_kT_{2,k}+\gamma_kT_{3,k}|\cF^k]=(1-\beta_k-\gamma_k)T_{1,k},
\end{align}
where we have used the condition $\mE[g_{\phi}(\bx)]=g(\bx)$ in the Assumption~\ref{assump:unbiased}.

Therefore, conditioned on $\cF^k$, taking expectation on the both sides of~\eqref{lemmaproof:1}, we have 
\begin{align*}
	&\mE[\norm{\by^{k+1}-g(\bx^k)}^2|\cF^k]\\
	=&\mE[\norm{(1-\beta_k)(\by^k-g(\bx^{k-1}))}^2|\cF^k]+\mE[\norm{(1-\beta_k)T_{1,k}+\beta_kT_{2,k}+\gamma_kT_{3,k}}^2|\cF^k]\\
	&+2\mE[\fprod{(1-\beta_k)(\by^k-g(\bx^{k-1})),(1-\beta_k)T_{1,k}+\beta_kT_{2,k}+\gamma_kT_{3,k}}|\cF^k]\\
	\leq &(1-\beta_k)^2\norm{\by^k-g(\bx^{k-1})}^2+2\mE[\norm{(1-\beta_k)T_{1,k}+\beta_kT_{2,k}}^2|\cF^k]+2\gamma_k^2\mE[\norm{T_{3,k}}^2|\cF^k]\\
	&+2(1-\beta_k)(1-\beta_k-\gamma_k)\norm{\by^{k}-g(\bx^{k-1})}\cdot\norm{T_{1,k}}\\
	\labelrel\leq{lemmaeq:1} &(1-\beta_k)^2\norm{\by^k-g(\bx^{k-1})}^2+2(1-\beta_k)^2\norm{T_{1,k}}^2+2\beta_k^2\mE[\norm{T_{2,k}}^2|\cF^k]+2\gamma_k^2\mE[|\norm{T_{3,k}}^2|\cF^K]\\
	&+(1-\beta_k)^2\beta_k\norm{\by^k-g(\bx^{k-1})}^2+\frac{(1-\beta_k-\gamma_k)^2}{\beta_k}\norm{T_{1,k}}^2\\
	\labelrel\leq{lemmaeq:2} & (1-\beta_k)\norm{\by^k-g(\bx^{k-1})}^2+2\beta_k^2V_g^2+(2(1-\beta_k)^2+2\gamma_k^2+\frac{(1-\beta_k-\gamma_k)^2}{\beta_k})\alpha_k^2C_g^4C_f^2
\end{align*}
where~\eqref{lemmaeq:1} is based on the Young's inequality;~\eqref{lemmaeq:2} uses the updating rule of $\by^{k+1}$ (see the Algorithm~\ref{alg:SCSC1}) and Assumptions~\ref{assump:unbiased} and~\ref{assumption:bounded1}.
\end{proof}

\subsubsection{Proof of Theorem~\ref{mainthm1}}
\begin{proof}
	Using the smoothness of $F(\bx^k)$ (see Assumption~\ref{assump:smooth_obj}), we have 
	\begin{align}
		&F(\bx^{k+1})\leq F(\bx^k)-\alpha_k\fprod{\eta^{k+1},\text{grad }F(\bx^k)}+\frac{\alpha_k^2 L_F}{2}\norm{\eta^{k+1}}^2\\
		=&F(\bx^k)-\alpha_k\norm{\text{grad }F(\bx^k)}^2+\frac{\alpha_k^2L_F}{2}\norm{\eta^{k+1}}^2\\
		+&\alpha_k\fprod{\text{grad} F(\bx^k),(Dg(\bx^k))^*\text{Proj}_{g(\bx^k)}\nabla f(g(\bx^k))-(Dg_{\phi^k}(\bx^k))^*\text{Proj}_{g_{\phi^k}(\bx^k)}\nabla f_{\xi^k}(\by^{k+1})},
	\end{align}
	where the first inequality  comes from the Proposition 10.47 in~\cite{boumal2020introduction}.
	
	Conditioned on $\cF^k$, taking expectation over $\phi^k$ and $\xi^k$ on both sides, we have 
	\begin{align}
		&\mE[F(\bx^{k+1})|\cF^k]\\
		\labelrel\leq{thmeq:1}&F(\bx^k)-\alpha_k\norm{\text{grad }F(\bx^k)}^2+\frac{\alpha_k^2L_F}{2}\mE[\norm{\eta^{k+1}}^2|\cF^k]\\
		+&\alpha_k \mE[\fprod{\text{grad} F(\bx^k),(Dg_{\phi^k}(\bx^k))^*\text{Proj}_{g_{\phi^k}(\bx^k)}(\nabla f_{\xi^k}(g(\bx^k))-\nabla f_{\xi^k}(\by^{k+1}))}|\cF^k]\\
		\labelrel\leq{thmeq:2}&F(\bx^k)-\alpha_k\norm{\text{grad }F(\bx^k)}^2+\frac{\alpha_k^2L_F}{2}\mE[\norm{\eta^{k+1}}^2|\cF^k]\\
		+&\alpha_k\norm{\text{grad} F(\bx^k)}\cdot\mE[\norm{Dg(\bx^k)[\cdot]}_{op}^2|\cF^k]^{1/2}\cdot\mE[\norm{\nabla f(g(\bx^k))-\nabla f(\by^{k+1})}^2|\cF^k]^{1/2}\\
		\labelrel\leq{thmeq:3}&F(\bx^k)-\alpha_k\norm{\text{grad }F(\bx^k)}^2+\frac{L_F}{2}\alpha_k^2C_g^2C_f^2\\
		+&\alpha_kC_g L_f\cdot\norm{\text{grad} F(\bx^k)}\cdot\mE[\norm{g(\bx^k)-\by^{k+1}}^2|\cF^k]^{1/2}\\
		\labelrel\leq{thmeq:4}&F(\bx^k)-\alpha_k(1-\frac{\alpha_k}{4\beta_k}C_g^2L_f^2)\norm{\text{grad }F(\bx^k)}^2+\frac{L_F}{2}\alpha_k^2C_g^2C_f^2+\beta_k\mE[\norm{g(\bx^k)-\by^{k+1}}^2|\cF^k].
	\end{align}
	The inequality~\eqref{thmeq:1} follows from Assumption~\ref{assump:oracle}; The inequality~\eqref{thmeq:2} uses the Cauchy-Schwartz inequality and nonexpansiveness property of the projection operator; The inequality~\eqref{thmeq:3} is based on Assumptions~\ref{assumption:smooth_f},~\ref{assumption:bounded1} and the Jensen inequality; Finally, the inequality~\eqref{thmeq:4} uses the Young's inequality. 	
	Based on the definition of the Lyapunov function in~\eqref{lyapunov}, it follows that 
	\begin{align}
		\mE[\cV^{k+1}|\cF^k]-\cV^k&\leq -\alpha_k(1-\frac{\alpha_k}{4\beta_k}C_g^2L_f^2)\norm{\text{grad }F(\bx^k)}^2+\frac{L_F}{2}\alpha_k^2C_g^2C_f^2\\
		&+(1+\beta_k)\mE[\norm{g(\bx^k)-\by^{k+1}}^2|\cF^k]-\norm{g(\bx^{k-1})-\by^k}^2.
	\end{align}
	In the following, we enforce $\beta_k=\frac{\alpha_kC_g^2L^2_f}{2}$, and $\alpha_k$ sufficiently small such that $\beta_k\in (0,1)$. 
	Combining the result in the Lemma~\ref{lemma:variance_tracking}, we have
	\begin{align}
		\mE[\cV^{k+1}|\cF^k]-\cV^k&\leq -\frac{\alpha_k}{2}\norm{\text{grad }F(\bx^k)}^2+\frac{L_F}{2}\alpha_k^2C_g^2C_f^2-\beta_k^2\norm{g(\bx^{k-1})-\by^k}^2\\
		&+\frac{1+\beta_k}{2}\alpha_k^2C_g^4L_f^4V_g^2+[2+6(t+1)^2]\alpha_k^2C_g^4C_f^2\\
		&\leq -\frac{\alpha_k}{2}\norm{\text{grad }F(\bx^k)}^2+\frac{L_F}{2}\alpha_k^2C_g^2C_f^2\\&\alpha_k^2C_g^4L_f^4V_g^2+[2+6(t+1)^2]\alpha_k^2C_g^4C_f^2
	\end{align}
	Defining $B:=\frac{L_F}{2}C_g^2C_f^2+C_g^4L_f^4V_g^2+[2+6(t+1)^2]C_g^4C_f^2$, and taking expectation over $\cF^k$ on both sides of the above inequality, it follows that  
	\begin{equation}
		\mE[\cV^{k+1}]\leq\mE[\cV^k]-\frac{\alpha_k}{2}\mE[\norm{\text{grad }F(\bx^k)}^2]+B\alpha_k^2.
	\end{equation}
	By telescoping, we have 
	\begin{equation}
		\mE[\cV^{K}]\leq\cV^0-\frac{1}{2}\sum_{k=0}^{K-1}\alpha_k\mE[\norm{\text{grad }F(\bx^k)}^2]+B\sum_{k=0}^{K-1}\alpha_k^2.
	\end{equation}.
	Using the fact that $\mE[\cV^K]\geq 0$ and rearranging the terms, we have 
	\begin{equation}
		\frac{\sum_{k=0}^{K-1}\alpha_k\mE[\norm{\text{grad }F(\bx^k)}^2]}{\sum_{k=0}^{K-1}\alpha_k}\leq \frac{2\cV^0+2B\sum_{k=0}^{K-1}\alpha_k^2}{\sum_{k=0}^{K-1}\alpha_k}.
	\end{equation}
	Choosing the stepsize as $\alpha_k=\frac{1}{\sqrt{{K}}}$ leads
	\begin{equation}
		\frac{\sum_{k=0}^{K-1}\mE[\norm{\text{grad }F(\bx^k)}^2]}{K}\leq \frac{2\cV^0+2B}{\sqrt{K}},
	\end{equation}
	from which the proof is complete. 
\end{proof}

\subsection{Complexity analysis of the multi-level R-SCGD} \label{sec:supp_multi-level}

\subsubsection{Proof of Lemma~\ref{lemma:variance_tracking_multi}}
\begin{proof}
	Let $\cF^{k,n}$ denotes the $\sigma$-algebra generated by $\cF^k\cup\{\theta^k_i,1\leq i\leq n-1\}$. From the update scheme for $\by_n^{k+1}$ (similarly when $n=1$), we have 
	\begin{align}
		\by_n^{k+1}-f_n(\by^{k+1}_{n-1})&=(1-\beta_k)(\by_n^k-f_n(\by^k_{n-1}))+(1-\beta_k)(f_n(\by^k_{n-1})-f_{n}(\by^{k+1}_{n-1}))\\
		&+\beta_k(f_{n}(\by^{k+1}_{n-1};\theta^k_n)-f_n(\by_{n-1}^{k+1}))+\gamma_k(f_n(\by^{k+1}_{n-1};\theta^k_n)-f_n(\by_{n-1}^{k};\theta^k_n))\\
		&=(1-\beta_k)(\by_n^k-f_n(\by_{n-1}^k))+(1-\beta_k)T_1+\beta_kT_2+\gamma_kT_3\label{lemma:multi1}
	\end{align}
	where $T_1:=f_n(\by^k_{n-1})-f_{n}(\by^{k+1}_{n-1})$, $T_2:=f_{n}(\by^{k+1}_{n-1};\theta^k_n)-f_n(\by_{n-1}^{k+1})$ and $T_3:=f_n(\by^{k+1}_{n-1};\theta^k_n)-f_n(\by_{n-1}^{k};\theta^k_n)$.
	
	Conditioned on $\cF^{k,n}$, taking expectation over $\theta^k_n$, we have 
	\begin{equation}\label{lemma:multi_equalzero}
		\mE[(1-\beta_k)T_1+\beta_2T_2+\gamma_k T_3|\cF^{k,n}]=(1-\beta_k-\gamma_k)T_1.
	\end{equation}
	where we have used Assumption~\ref{assump:unbiased_multi}. 
	
	Based on the~\eqref{lemma:multi1}, we have 
	\begin{align*}
		&\mE[\norm{\by_n^{k+1}-f_n(\by_{n-1}^{k+1})}^2|\cF^{k,n}]\\
		\labelrel={lemmamultieq:1}&(1-\beta_k)^2\norm{\by^k_n-f_n(\by^k_{n-1})}^2+\mE[\norm{(1-\beta_k)T_1+\beta_kT_2+\gamma_kT_3}^2|\cF^{k,n}]\\
		&+2(1-\beta_k)(1-\beta_k-\gamma_k)\norm{\by^k_n-f_n(\by^k_{n-1})}\cdot\norm{T_1}\\
		\leq&(1-\beta_k)^2\norm{\by^k_n-f_n(\by^k_{n-1})}^2+2\mE[\norm{(1-\beta_k)T_1+\beta_kT_2}^2|\cF^{k,n}]+2
		\gamma_k^2\mE[\norm{T_3}^2|\cF^{k,n}]\\
		&+(1-\beta_k)^2\beta_k\norm{\by^k_n-f_n(\by^k_{n-1})}^2+\frac{(1-\beta_k-\gamma_k)^2}{\beta_k}\norm{T_1}^2\\
		\labelrel\leq{lemmamultieq:2}&(1-\beta_k)\norm{\by^k_n-f_n(\by^k_{n-1})}^2+[2(1-\beta_k)^2+\frac{(1-\beta_k-\gamma_k)^2}{\beta_k}] \norm{T_1}^2++2\gamma_k^2\mE[\norm{T_3}^2|\cF^{k,n}]\\
		+&2\beta_k^2\mE[\norm{T_2}|\cF^{k,n}]\\
		\labelrel\leq{lemmamultieq:3}&(1-\beta_k)\norm{\by^k_n-f_n(\by^k_{n-1})}^2+[2(1-\beta_k)^2+\frac{(1-\beta_k-\gamma_k)^2}{\beta_k}+2\gamma_k^2]C_n^2\norm{\by^k_{n-1}-\by^{k+1}_{n-1}}^2+2\beta_k^2V_n^2,
	\end{align*}
	where~\eqref{lemmamultieq:1} follows from the equality~\eqref{lemma:multi_equalzero};~\eqref{lemmamultieq:2} is based on the fact that $$\mE[\fprod{T_1,T_2}|\cF^{k,n}]=\fprod{T_1,\mE[T_2|\cF^{k,n}]}=0\text{ and }(1-\beta_k)^2(1+\beta_k)\leq 1-\beta_k;$$
	and the last inequality~\eqref{lemmamultieq:3} follows from Assumptions~\ref{assump:unbiased_multi} and~\ref{assumption:bounded1_multi}.
	
Due to the law of total expectation (see e.g.~\cite{billingsley2008probability}), 
	\begin{align*}
		&\mE[\norm{\by_n^{k+1}-f_n(\by_{n-1}^{k+1})}^2|\cF^k]\\
		=&\mE[\mE[\norm{\by_n^{k+1}-f_n(\by_{n-1}^{k+1})}^2|\cF^{k,n}]|\cF^k]\\
		\leq&(1-\beta_k)\norm{\by^k_n-f_n(\by^k_{n-1})}^2+[2(1-\beta_k)^2+\frac{(1-\beta_k-\gamma_k)^2}{\beta_k}+2\gamma_k^2]C_n^2\norm{\by^k_{n-1}-\by^{k+1}_{n-1}}^2+2\beta_k^2V_n^2.
	\end{align*}
\end{proof}

\subsubsection{Proof of Theorem~\ref{mainthm2}}
\begin{proof}
	Using the smoothness of $F(\bx^k)$ (see Assumption~\ref{assump:smooth_obj}), we have 
	\begin{align*}
		F(\bx^{k+1})&\leq F(\bx^k)-\alpha_k\fprod{\eta^{k+1},\text{grad }F(\bx^k)}+\frac{\alpha_k^2 L_F}{2}\norm{\eta^{k+1}}^2\\
		&=F(\bx^k)-\alpha_k\norm{\text{grad }F(\bx^k)}^2+\frac{\alpha_k^2L_F}{2}\norm{\eta^{k+1}}^2
		+\alpha_k\fprod{\text{grad} F(\bx^k),\text{grad} F(\bx^k)-\eta^{k+1}},
	\end{align*}
	where the first inequality  comes from the Proposition 10.47 in~\cite{boumal2020introduction}.
	Conditioned on $\cF^k$, taking expectation over $\theta_1^k$,...,$\theta_N^k$, we have 
	\begin{align*}
		&\mE[F(\bx^{k+1})|\cF^k]\\
		\leq & F(\bx^k)-\alpha_k\norm{\text{grad} F(\bx^k)}^2+\frac{\alpha_k^2 L_F}{2}\mE[\norm{\eta^{k+1}}^2|\cF^k]+\alpha_k\fprod{\text{grad} F(\bx^k),\text{grad} F(\bx^k)-\mE[\eta^{k+1}|\cF^k]}\\
		\leq &F(\bx^k)-\alpha_k\norm{\text{grad} F(\bx^k)}^2+\frac{\alpha_k^2L_F}{2}\prod_{i=n}^N C_n^2+\alpha_k\norm{\text{grad} F(\bx^k)}\cdot\norm{\text{grad} F(\bx^k)-\mE[\eta^{k+1}|\cF^k]}\\
		\labelrel\leq{mainthm2:2}& F(\bx^k)-\alpha_k\norm{\text{grad} F(\bx^k)}^2+\frac{\alpha_k^2L_F}{2}\prod_{n=1}^N C_n^2+\alpha_k\sum_{n=1}^{N-1}A_n\norm{\text{grad} F(\bx^k)}\cdot \mE[\norm{\by_n^{k+1}-f_n(\by^{k+1}_{n-1})}|\cF^k]\\
		\labelrel\leq{mainthm2:3}& F(\bx^k)-\alpha_k(1-\frac{\alpha_k}{4\beta_k}\sum_{n=1}^{N-1}A_n^2)\norm{\text{grad} F(\bx^k)}^2+\beta_k\sum_{n=1}^{N-1}\mE[\norm{\by_n^{k+1}-f_n(\by_{n-1}^{k+1})}^2|\cF^k]+\frac{\alpha_k^2L_F}{2}\prod_{n=1}^{N}C_n^2, 
	\end{align*}
	where $A_n:=\sum_{m=n+1}^{N-1}C_N\cdots C_{m+1}C_{m-1}\cdots C_1L_m\cdots L_{n+1}$. Inequality~\eqref{mainthm2:2} is based on Assumption~\ref{assump:oracle_multi} and follows exactly the proof in~\cite{chen2021solving}; the inequality~\eqref{mainthm2:3} follows from the Young's inequality.
	Based on the definition of generalized Lyapunov function~\eqref{lyapunov:multi}, we have 
	\begin{align*}
		\mE[\cV^{k+1}|\cF^k]&-\cV^k\leq \mE[F(\bx^{k+1})|\cF^k]+\mE[\sum_{n=1}^{N-1}\norm{\by_n^{k+1}-f_n(\by_{n-1}^{k+1})}^2|\cF^k]-F(\bx^k)-\sum_{n=1}^{N-1}\norm{\by_n^k-f_n(\by_{n-1}^k)}^2 \\
		\leq &-\alpha_k(1-\frac{\alpha_k}{4\beta_k}\sum_{n=1}^{N-1}A_n^2)\norm{\text{grad} F(\bx^k)}^2+(\beta_k+1)\sum_{n=1}^{N-1}\mE[\norm{\by_n^{k+1}-f_n(\by_{n-1}^{k+1})}^2|\cF^k]\\
		&+\frac{\alpha_k^2L_F}{2}\prod_{n=1}^{N}C_n^2-\sum_{n=1}^{N-1}\norm{\by_n^k-f_n(\by_{n-1}^k)}^2\\
		&=-\frac{\alpha_k}{2}\norm{\text{grad} F(\bx^k)}^2+(2\beta_k+1)\sum_{n=1}^{N-1}\mE[\norm{\by_n^{k+1}-f_n(\by_{n-1}^{k+1})}^2|\cF^k]+\frac{\alpha_k^2L_F}{2}\prod_{n=1}^{N}C_n^2\\
		&-\beta_k\sum_{n=1}^{N-1}\mE[\norm{\by_n^{k+1}-f_n(\by_{n-1}^{k+1})}^2|\cF^k]-\sum_{n=1}^{N-1}\norm{\by_n^k-f_n(\by_{n-1}^k)}^2, 
	\end{align*}
	where in the last equality we set
	\begin{equation}\label{constantchoise1}
		\beta_k:=\frac{\alpha_k}{2}\sum_{n=1}^{N-1}A_n^2.
	\end{equation}
	Before we combine the above inequality with the result in Lemma~\ref{lemma:variance_tracking_multi}, we consider the following quantities, 
	\begin{equation}\label{quantities_main}
		\begin{split}
			&(2\beta_k+1)(1-\beta_k)\leq 1 \text{ when } 0\leq\beta_k\leq 1,\\
			&[2(1-\beta_k)^2+\frac{(1-\beta_k-\gamma_k)^2}{\beta_k}+2\gamma_k^2] (2\beta_k+1)\leq M_1:= 2+3(t+1)(t+2),
		\end{split}
	\end{equation}
	where we choose 
	\begin{equation}\label{constantchoise2}
		\gamma_k=1-t_k\beta_k \text{ such that }t:=\sup\{|t_k|\} \text{ is finite}.
	\end{equation}
	Using the quantities~\eqref{quantities_main} and Lemma~\ref{lemma:variance_tracking_multi}, we have 
	\begin{equation}\label{mainthm2:eq9}
	\begin{split}
		&\mE[\cV^{k+1}|\cF^k]-\cV^k \\
		&\leq -\frac{\alpha_k}{2} \norm{\text{grad} F(\bx^k)}^2+M_1\sum_{n=1}^{N-1}C_n^2\mE[\norm{\by^k_{n-1}-\by^{k+1}_{n-1}}^2|\cF^k]+2\beta_k^2(2\beta_k+1)\sum_{n=1}^{N-1}V_n^2\\
		&+\frac{\alpha_k^2L_F}{2}\prod_{n=1}^N C_n^2-\beta_k\sum_{n=1}^{N-1}\mE[\norm{\by_n^{k+1}-f_n(\by_{n-1}^{k+1})}^2|\cF^k]\\
		&=-\frac{\alpha_k}{2} \norm{\text{grad} F(\bx^k)}^2+M_1 C_1^2\alpha_k^2\mE[\norm{\eta^k}^2|\cF^k]+2\beta_k^2(2\beta_k+1)\sum_{n=1}^{N-1}V_n^2+\frac{\alpha_k^2 L_F}{2}\prod_{n=1}^N C_n^2 \\
		&+\sum_{n=2}^{N-1}(M_1 C_n^2+\psi_n)\mE[\norm{\by_{n-1}^k-\by_{n-1}^{k+1}}^2|\cF^k]
		-\sum_{n=2}^{N-1}\psi_n\mE[\norm{\by_{n-1}^k-\by_{n-1}^{k+1}}^2|\cF^k] \\
		&-\beta_k\sum_{n=1}^{N-1}\mE[\norm{\by_n^{k+1}-f_n(\by_{n-1}^{k+1})}^2|\cF^k].
	\end{split}
	\end{equation}
	The parameters $\psi_n$ will be carefully chosen such that the last three terms are upper bounded by certain constants. 
	Due to the updating rule for $\by_n^{k+1}$ in Algorithm~\ref{alg:SCSC2}, we have 
	\begin{equation}
		(1-\beta_k)(\by_{n-1}^{k+1}-\by_{n-1}^k)=\beta_k(f_{n-1}(\by_{n-2}^{k+1};\theta_{n-1}^k)-\by_{n-1}^{k+1})+\gamma_k(f_{n-1}(\by_{n-2}^{k+1};\theta_{n-1}^k)-f_{n-1}(\by_{n-2}^k;\theta_{n-1}^k))
	\end{equation}
	Conditioned on $\cF^k$, we have 
	\begin{equation}\label{mainthm2:eq10}
	\begin{split}
		&\mE[\norm{\by_{n-1}^{k+1}-\by_{n-1}^k}^2|\cF^k]\\
		&\leq 2(\frac{\beta_k}{1-\beta_k})^2\mE[\norm{f_{n-1}(\by_{n-2}^{k+1};\theta_{n-1}^k)-f_{n-1}(\by_{n-2}^{k+1})+f_{n-1}(\by_{n-2}^{k+1})-\by_{n-1}^{k+1}}^2|\cF^k]\\
		& +2(\frac{\gamma_k}{1-\beta_k})^2\mE[\norm{f_{n-1}(\by_{n-2}^{k+1};\theta_{n-1}^k)-f_{n-1}(\by_{n-2}^k;\theta_{n-1}^k)}^2|\cF^k]\\
		&\leq 2(\frac{\beta_k}{1-\beta_k})^2V_{n-1}^2+2(\frac{\beta_k}{1-\beta_k})^2\mE[\norm{f_{n-1}(\by_{n-2}^{k+1})-\by_{n-1}^{k+1}}^2|\cF^k]\\&+2(\frac{\gamma_k}{1-\beta_k})^2C_{n-1}^2\mE[\norm{\by_{n-2}^{k+1}-\by_{n-2}^k}^2|\cF^k]\\
		&\leq 2(\frac{\beta_k}{1-\beta_k})^2V_{n-1}^2+2(\frac{\beta_k}{1-\beta_k})^2\mE[\norm{f_{n-1}(\by_{n-2}^{k+1})-\by_{n-1}^{k+1}}^2|\cF^k]\\&+2t^2C_{n-1}^2\mE[\norm{\by_{n-2}^{k+1}-\by_{n-2}^k}^2|\cF^k].
	\end{split}
	\end{equation}
	Plugging~\eqref{mainthm2:eq10} into~\eqref{mainthm2:eq9}, we have 
	\begin{equation}
	\begin{split}
		\mE[\cV^{k+1}|\cF^k]&-\cV^k\leq-\frac{\alpha_k}{2} \norm{\text{grad} F(\bx^k)}^2+M_1 C_1^2\alpha_k^2\mE[\norm{\eta^k}^2|\cF^k]+2\beta_k^2(2\beta_k+1)\sum_{n=1}^{N-1}V_n^2 \\
		&+\frac{\alpha_k^2 L_F}{2}\prod_{n=1}^N C_n^2
		-\sum_{n=2}^{N-1}\psi_n\mE[\norm{\by_{n-1}^k-\by_{n-1}^{k+1}}^2|\cF^k]-\beta_k\sum_{n=1}^{N-1}\mE[\norm{\by_n^{k+1}-f_n(\by_{n-1}^{k+1})}^2|\cF^k]\\
		&+2(\frac{\beta_k}{1-\beta_k})^2 \sum_{n=2}^{N-1}(M_1 C_n^2+\psi_n)V_{n-1}^2 \\
		&+2(\frac{\beta_k}{1-\beta_k})^2\sum_{n=2}^{N-1}(M_1 C_n^2+\psi_n)\mE[\norm{f_{n-1}(\by_{n-2}^{k+1})-\by_{n-1}^{k+1}}^2|\cF^k]\\&+2t^2\sum_{n=2}^{N-1}(M_1 C_n^2+\psi_n)C_{n-1}^2\mE[\norm{\by_{n-2}^{k+1}-\by_{n-2}^k}^2|\cF^k]
	\end{split}
	\end{equation}
	Choosing parameters $\{\psi_n\}$ and $\{\beta_k\}$ such that 
	\begin{equation}\label{constantchoise3}
	\begin{split}
	2t^2(M_1C_n^2+\psi_n)C_{n-1}^2&\leq \psi_{n-1},\\
	2(\frac{\beta_k}{1-\beta_k})^2(M_1C_n^2+\psi_n)&\leq \beta_k,\\
	\beta_k\leq 1/2,
	\end{split}
	\end{equation}
	
	(for instance,~\eqref{constantchoise3} can be satisfied by choosing $\psi_{N-1}=0$, $\psi_{N-2}=2t^2M_1C_{N-1}^2C_{N-2}^2$, ... 
	and $\beta_k$ sufficiently small),
	we have 
	\begin{equation}
		\begin{split}
		\mE[\cV^{k+1}|\cF^k]&-\cV^k\leq-\frac{\alpha_k}{2} \norm{\text{grad} F(\bx^k)}^2+M_1 C_1^2\alpha_k^2\mE[\norm{\eta^k}^2|\cF^k]+2\beta_k^2(2\beta_k+1)\sum_{n=1}^{N-1}V_n^2 \\
		&+\frac{\alpha_k^2 L_F}{2}\prod_{n=1}^N C_n^2
		+2(\frac{\beta_k}{1-\beta_k})^2 \sum_{n=2}^{N-1}(M_1 C_n^2+\psi_n)V_{n-1}^2\\
		&\leq-\frac{\alpha_k}{2} \norm{\text{grad} F(\bx^k)}^2+ M_1C_1^2\alpha_k^2 \prod_{n=1}^N C_n^2+\frac{\alpha_k^2 L_F}{2}\prod_{n=1}^N C_n^2+2\beta_k^2(2\beta_k+1)\sum_{n=1}^{N-1}V_n^2\\
		&+2(\frac{\beta_k}{1-\beta_k})^2 \sum_{n=2}^{N-1}(M_1 C_n^2+\psi_n)V_{n-1}^2\\
		&\leq -\frac{\alpha_k}{2} \norm{\text{grad} F(\bx^k)}^2+ M_2 \alpha_k^2 +M_3 \beta_k^2,
		\end{split}
	\end{equation}
	where $M_2:=M_1C_1^2\prod_{n=1}^NC_n^2+\frac{L_F}{2}\prod_{n=1}^NC_n^2$ and $M_3:=4\sum_{n=1}^{N-1}V_n^2+8\sum_{n=2}^{N-1}(M_1C_n^2+\psi_n)V_{n-1}^2$. 
	Choosing the stepsize 
		$\alpha_k=1/\sqrt{K}$,
we have 
\begin{equation}
	\frac{\sum_{k=0}^{K-1}\mE[\norm{\text{grad} F(\bx^k)}^2]}{K}\leq \frac{2\cV^0+2(M_2+M_3(\sum_{n=1}^{N-1}A_n^2)^2/4)}{\sqrt{K}},
\end{equation}
given that the parameters are chosen as~\eqref{constantchoise1},\eqref{constantchoise2},\eqref{constantchoise3}.
	
\end{proof}

\end{document}